\documentclass[10pt]{article}

\usepackage{amssymb}
\usepackage[margin=1.25in]{geometry}
\usepackage{tikz,amsmath,listings,textcomp}
\usepackage{xfrac,subfigure}

\newcommand{\e}{\ensuremath{\varepsilon}}

\title{Backward Error Analysis for Perturbation Methods
\thanks{We would like to thank Pei Yu, Robert Moir and Julia Jankowski for their various contributions to this paper. We are also indebted to NSERC, Western University, as well as Galima Hassan for the key logistic support they provided.} 
}

\author{
  Corless, Robert M.\\
  Western University\\
  \texttt{rcorless@uwo.ca}
  \and
  Fillion, Nicolas\\
  Simon Fraser University\\
  \texttt{nfillion@sfu.ca}
}
\date{\today}

\begin{document}

\maketitle

\begin{abstract}
We demonstrate via several examples how the backward error viewpoint can be used in the analysis of solutions obtained by perturbation methods. We show that this viewpoint is quite general and offers several important advantages. Perhaps the most important is that backward error analysis can be used to demonstrate the validity of the solution, however obtained and by whichever method. This includes a nontrivial safeguard against slips, blunders, or bugs in the original computation. We also demonstrate its utility in deciding when to truncate an asymptotic series, improving on the well-known rule of thumb indicating truncation just prior to the smallest term. We also give an example of elimination of \textsl{spurious} secular terms even when genuine secularity is present in the equation. We give short expositions of several well-known perturbation methods together with computer implementations (as scripts that can be modified). We also give a generic backward error based method that is equivalent to iteration (but we believe useful as an organizational viewpoint) for regular perturbation. 
\end{abstract}


\section{Introduction}
As the title suggests, the main idea of this paper is to use backward error analysis (BEA) to assess and interpret solutions obtained by perturbation methods. The idea will seem natural, perhaps even obvious, to those who are familiar with the way in which backward error analysis has seen its scope increase dramatically since the pioneering work of Wilkinson in the 60s, e.g., \cite{Wilkinson(1963),Wilkinson(1965)}.   From early results in numerical linear algebraic problems and computer arithmetic, it has become a general method fruitfully applied to problems involving root finding, interpolation,  numerical differentiation, quadrature, and the numerical solutions of ODEs, BVPs, DDEs, and PDEs, see, e.g, \cite{CorlessFillion(2013),Deuflhard(2003),Higham(1996)}. This is hardly a surprise when one considers that BEA offers several interesting advantages over a purely forward-error approach.

BEA is often used in conjunction with perturbation methods. Not only is it the case that many algorithms' backward error analyses rely on perturbation methods, but the backward error is related to the forward error by a coefficient of sensitivity known as the condition number, which is itself a kind of sensitivity to perturbation. In this paper, we examine an apparently new idea, namely, that perturbation methods themselves can also be interpreted within the backward error analysis framework. Our examples will have a classical feel, but the analysis and interpretation is what differs, and we will make general remarks about the benefits of this mode of analysis and interpretation. 

However, due to the breadth of the literature in perturbation theory, we cannot determine with certainty the extent to which applying backward error analysis to perturbation methods is new. Still, none of the works we know, apart from \cite{Boyd(2014)}, \cite{Corless(1993)b}, and \cite{Corless(2014)}, even mention the possibility of using of using BEA to explain or measure the success of a perturbation computation. Among the books we have consulted, only \cite[p.~251 \& p.~289]{Boyd(2014)} mentions the residual by name, but does not use it systematically. At the very least, therefore, the idea of using BEA in relation to perturbation methods might benefit from a wider discussion.

\section{The basic method from the BEA point of view} \label{genframe}

The basic idea of BEA is increasingly well-known in the context of numerical methods. The slogan \textsl{a good numerical method gives the exact solution to a nearby problem} very nearly sums up the whole perspective. Any number of more formal definitions and discussions exist---we like the one given in \cite[chap.~1]{CorlessFillion(2013)}, as one might suppose is natural, but one could hardly do better than go straight to the source and consult, e.g., \cite{Wilkinson(1963),Wilkinson(1965),Wilkinson(1971),Wilkinson(1984)}. More recently \cite{Grcar(2011)} has offered a good historical perspective. In what follows we give a brief formal presentation and then give detailed analyses by examples in subsequent sections.


Problems can generally be represented as maps from an input space $\mathcal{I}$ to an output space $\mathcal{O}$.
  If we have a problem $\varphi:\mathcal{I}\to\mathcal{O}$ and wish to find $y=\varphi(x)$ for some putative input $x\in\mathcal{I}$, lack of tractability might instead lead you to engineer a simpler problem $\hat{\varphi}$ from which you would compute $\hat{y}=\hat{\varphi}(x)$. Then $\hat{y}-y$ is the \textsl{forward error} and, provided it is small enough for your application, you can treat $\hat{y}$ as an approximation in the sense that $\hat{y}\approx \varphi(x)$. In BEA, instead of focusing on the forward error, we try to find an $\hat{x}$ such that $\hat{y}=\varphi(\hat{x})$ by considering the \textsl{backward error} $\Delta x=\hat{x}-x$, i.e., we try to find for which set of data our approximation method $\hat{\varphi}$ has exactly solved our reference problem $\varphi$. The general picture can be represented by the following commutative diagram:
  \begin{center}
  \begin{tikzpicture}
  \def\a{2}
  \draw (0,0) node (x) {$x$};
  \draw (\a,0) node (y) {$y$};
  \draw (0,-\a) node (xhat) {$\hat{x}$};
  \draw (\a,-\a) node (yhat) {$\hat{y}$};
  \draw (x) edge[->] node[above] {$\varphi$} (y);
  \draw (xhat) edge[->,dashed] node[below] {$\varphi$} (yhat);
  \draw (x) edge[->,dashed] node[left] {$+\Delta x$} (xhat);
\draw (y) edge[->] node[right] {$+\Delta y$} (yhat);
\draw (x) edge[->] node[above right] {$\hat{\varphi}$} (yhat);
  \end{tikzpicture}
  \end{center}
 We can see that, whenever $x$ itself has many components, different backward error analyses will be possible since we will have the option of reflecting the forward error back into different selections of the components. 
  
It is often the case that the map $\varphi$ can be defined as the solution to $\phi(x,y)=0$ for some operator $\phi$, i.e., as having the form
\begin{align} x\xrightarrow{\varphi} \left\{ y\mid \phi(x,y)=0\right\}\>.\end{align}
In this case, there will in particular be a simple and useful backward error resulting from computing the residual $r=\phi(x,\hat{y})$. Trivially $\hat{y}$ then exactly solves the reverse-engineered problem $\hat{\varphi}$ given by $\hat{\phi}(x,y)=\phi(x,y)-r=0$.          
Thus, when the residual can be used as a backward error, this directly computes a reverse-engineered problem that our method has solved exactly. We are then in the fortunate position of having both a problem and its solution, and the challenge then consists in determining how similar the reference problem $\varphi$ and the modified problems $\hat{\varphi}$ are, \textsl{and whether or not the modified problem is a good model for the phenomenon being studied}.

\paragraph{Regular perturbation BEA-style}
Now let us introduce a \textsl{general framework for perturbation methods} that relies on the general framework for BEA introduced above.
%
Perturbation methods are so numerous and varied, and the problems tackled are from so many areas, that it seems a general scheme of solution would necessarily be so abstract as to be difficult to use in any particular case. 
Actually, the following framework covers many methods. For simplicity of exposition, we will introduce it using the simple gauge functions $1,\e,\e^2,\ldots$, but note that extension to other gauges is usually straightforward (such as Puiseux, $\e^n\ln^m\e$, etc), as we will show in the examples. 
To begin with, let
\begin{align}
F(x,u;\e)=0 \label{operatoreq}
\end{align}
be the operator equation we are attempting to solve for the unknown $u$. The dependence of $F$ on the scalar parameter $\e$ and on any data $x$ is assumed but henceforth not written explicitly. In the case of a simple power series perturbation, we will take the $m$th order approximation to $u$ to be given by the \textsl{finite} sum
\begin{align}
z_m = \sum_{k=0}^m \e^ku_k\>.
\end{align}
The operator $F$ is assumed to be Fr\'echet differentiable. For convenience we assume slightly more, namely, that for any $u$ and $v$ in a suitable region, there exists a linear invertible operator $F_1(v)$ such that
\begin{align}
F(u) = F(v) + F_1(v)(u-v) + O\left(\|u-v\|^2\right)\>.
\end{align}
Here, $\|\cdot\|$ denotes any convenient norm. We denote the \textsl{residual} of $z_m$ by
\begin{align}
\Delta_m := F(z_m)\>,
\end{align}
\emph{i.e.}, $\Delta_m$ results from evaluating $F$ at $z_m$ instead of evaluating it at the reference solution $u$ as in equation \eqref{operatoreq}. If $\|\Delta_m\|$ is small, we say we have solved a ``nearby'' problem, namely, the reverse-engineered problem for the unknown $u$ defined by
\begin{align}
F(u)-F(z_m) = 0\>,
\end{align}
which is exactly solved by $u=z_m$. Of course this is trivial. It is \textsl{not} trivial in consequences if $\|\Delta_m\|$ is small compared to data errors or modelling errors in the operator $F$. We will exemplify this point more concretely later.

We now suppose that we have somehow found $z_0=u_0$, a solution with a residual whose size is such that
\begin{align}
\|\Delta_0\|=\|F(u_0)\| = O(\e)\qquad \textrm{as} \qquad \e\to0\>.
\end{align}
Finding this $u_0$ is part of the art of perturbation; much of the rest is mechanical.
Suppose now inductively that we have found $z_n$ with residual of size
\[
\|\Delta_n\| = O\left(\e^{n+1}\right) \quad\textrm{ as }\quad \e\to0\>.
\]
Consider $F(z_{n+1})$ which, by definition, is just $F(z_n+\e^{n+1}u_{n+1})$. We wish to choose the term $u_{n+1}$ in such a way that $z_{n+1}$ has residual of size $\|\Delta_{n+1}\|=O(\e^{n+2})$ as $\e\to0$. Using the Fr\'echet derivative of the residual of $z_{n+1}$ at $z_n$, we see that
\begin{align}
\Delta_{n+1} &= F(z_n+\e^{n+1}u_{n+1})= F(z_n)+F_1(z_n)\e^{n+1}u_{n+1}+O\left(\e^{2n+2}\right)\>. \label{resseries1}
\end{align}
By linearity of the Fr\'echet derivative, we also obtain $F_1(z_n) = F_1(z_0)+O(\e)= [\e^0]F_1(z_0)+O(\e)$. Here, $[\e^k]G$ refers to the coefficient of $\e^k$ in the expansion of $G$. Let
\begin{align}
A=[\e^0]F_1(z_0)\>,
\end{align}
 that is, the zeroth order term in $F_1(z_0)$. Thus, we reach the following expansion of $\Delta_{n+1}$:
\begin{align}
\Delta_{n+1} = F(z_n) + A\e^{n+1}u_{n+1}+O\left(\e^{n+2}\right)\>.\label{eqDnp1}
\end{align}
Note that, in equation \eqref{resseries1},  one could keep $F_1(z_n)$, not simplifying to $A$ and compute not just $u_{n+1}$ but, just as in Newton's method, double the number of correct terms. However, this in practice is often too expensive \cite[chap.~6]{Geddes(1992)b}, and so we will in general use this simplification. As noted, we only need $F_1(z_0)$ accurate to $O(\e)$, so in place of $F_1(z_0)$ in equation \eqref{eqDnp1} we use $A$. 

As a result of the above expansion of $\Delta_{n+1}$, we now see that to make $\Delta_{n+1} = O\left(\e^{n+2}\right)$, we must have $F(z_n)+A\e^{n+1}u_{n+1}=O(\e^{n+2})$, in which case
\begin{align}
A u_{n+1} +\frac{F(z_n)}{\e^{n+1}} = Au_{n+1} +\frac{\Delta_n}{\e^{n+1}} =O(\e)\>.
\end{align}
Since by hypothesis $\Delta_n=F(z_n)=O(\e^{n+1})$, we know that $\sfrac{\Delta_n}{\e^{n+1}}=O(1)$.
In other words, to find $u_{n+1}$ we solve the linear operator equation
\begin{align*}
A u_{n+1} = 
 -[\e^{n+1}]\Delta_n\>,
\end{align*}
where, again, $[\e^{n+1}]$ is the coefficient of the $(n+1)$th power of $\e$ in the series expansion of $\Delta$. Note that by the inductive hypothesis the right hand side has norm $O(1)$ as $\e\to0$. Then $\|\Delta_{n+1}\| = O(\e^{n+2})$ as desired, so $u_{n+1}$ is indeed the coefficient we were seeking.
We thus need $A=[\e^0]F(z_0)$ to be invertible. If not, the problem is singular, and essentially requires reformulation.\footnote{We remark that it is a sufficient but not necessary condition for regular expansion to be able to find our initial point $u_0$ and to have invertible $A=F_1(u_0;0)$. A regular perturbation problem can be defined in many ways, not just in the way we have done, with invertible $A$. For example, \cite[Sec 7.2]{Bender(1978)} essentially uses continuity in $\e$ as $\e\to0$ to characterize it. Another characterization is that for regular perturbation problems infinite perturbation series are convergent for some non-zero radius of convergence. 
} We shall see examples. If $A$ is invertible, the problem is regular.



This general scheme can be compared to that of, say, \cite{Bellman(1972)}. Essential similarities can be seen. In Bellman's treatment, however, the residual is used implicitly, but not named or noted, and instead the equation defining $u_{n+1}$ is derived by postulating an infinite expansion
\begin{align}
u=u_0+\e u_1+\e^2u_2+\cdots\>.
\end{align}
By taking the coefficient of $\e^{n+1}$ in the expansion of $\Delta_n$ we are implicitly doing the same work, but we will see advantages of this point of view. %
Also, note that in the frequent case of more general asymptotic sequences, namely Puiseux series or generalized approximations containing logarithmic terms, we can make the appropriate changes in a straightforward manner, as we will show below.

\section{Algebraic equations}

We begin by applying the regular method from section \ref{genframe} to algebraic equations. We begin with a simple scalar equation and gradually increase the difficulty, thereby demonstrating the flexibility of the backward error point of view.

\subsection{Regular perturbation}\label{RegularPert}
In this section, after applying the method from section \ref{genframe} to a scalar equation, we use the same method to solve a $2\times2$ system; higher dimensional systems can be solved similarly. We give some computer algebra implementations (scripts that the reader may modify) of the basic method. Finally, in this section, we give an alternative method based on the Davidenko equation that is simpler to use in Maple.

\subsubsection{Scalar equations}
Let us consider a simple example similar to many used in textbooks for classical perturbation analysis. Suppose we wish to find a real root of
\begin{align}
x^5 -x-1=0 \label{refprobalgeq}
\end{align}
and, since the Abel-Ruffini theorem---which says that in general there are no solutions in radicals to equations of degree 5 or more---suggests it is unlikely that we can find an elementary expression for the solution of this \textsl{particular} equation of degree 5, we introduce a parameter which we call $\e$, and moreover which we suppose to be small. That is, we embed our problem in a parametrized family of similar problems. If we decide to introduce $\e$ in the degree-1 term, so that 
\begin{align}
u^5-\e u-1=0\>, \label{pertalgeq}
\end{align}
we will see that we have a so-called regular perturbation problem. 

To begin with, we wish to find a $z_0$ such that $\Delta_0=F(z_0) = z_0^5-\e z_0-1=O(\e)$. Quite clearly, this can happen only if $z_0^5-1=0$. Ignoring the complex roots in this example, we take $z_0=1$. To continue the solution process, we now suppose that we have found
\begin{align}
z_n = \sum_{k=0}^n u_k\e^k
\end{align}
such that $\Delta_n=F(z_n) = z_n^5-\e z_n-1=O(\e^{n+1})$ and we wish to use our iterative procedure. We need the Fr\'echet derivative of $F$, which in this case is just
\begin{align}
F_1(u) &= 5u^4-\e\>,
\end{align}
because
\begin{align}
F(u) = u^5-\e u-1 &= v^5-\e v-1 + F'(v)(u-v)+O(u-v)^2\>.
\end{align}
Hence, $A=5z_0^4=5$, which is invertible. As a result our iteration is $\Delta_n=F(z_n)$, i.e., 
\begin{align}
5u_{n+1} = -[\e^{n+1}]\Delta_n\>.
\end{align}
Carrying out a few steps we have
\begin{align}
\Delta_0 = F(z_0) = F(1) = 1-\e-1 = -\e
\end{align}
so
\begin{align}
5\cdot u_1 = -[\e]\Delta_0 = -[\e](-\e) = 1\>.
\end{align}
Thus, $u_1=\sfrac{1}{5}$. Therefore, $z_1=1+\sfrac{\e}{5}$ and
\begin{align}
\Delta_1 &= \left(1+\frac{\e}{5}\right)^5 -\e\left(1+\frac{\e}{5}\right)-1\\
& = \left(1+5\frac{\e}{5}+10\frac{\e^2}{25}+O\left(\e^3\right)\right) - \e-\frac{\e^2}{5}-1\\
&=\left(\frac{2}{5}-\frac{1}{5}\right)\e^2+O\left(\e^3\right) = \frac{1}{5}\e^2+O\left(\e^3\right)\>.
\end{align}
Then we find that $Au_1=-\sfrac{1}{5}$ and thus $u_1=-\sfrac{1}{25}$. So, $u=1+\sfrac{\e}{5}-\sfrac{\e^2}{25}+O(\e^3)$. Finding more terms by this method is clearly possible although tedium might be expected at higher orders.
%
%
Luckily nowadays computers and programs are widely available that can solve such problems without much human effort, but before we demonstrate that, let's compute the residual of our computed solution so far:
\[
z_2 = 1+\frac{1}{5}\e-\frac{1}{25}\e^2\>.
\]
Then $\Delta_2 = z_2^5-\e z_2-1$ is
\begin{align}
\Delta_2  &= \left(1+\frac{1}{5}\e-\frac{1}{25}\e^2\right)^5-\e\left(1+\frac{1}{5}\e-\frac{1}{25}\e^2\right)-1 \nonumber \\
& = -\frac{1}{25}\e^3 - \frac{3}{125}\e^4+\frac{11}{3125}\e^5 +\frac{3}{125}\e^6 -\frac{2}{15625}\e^7 \nonumber\\ &\qquad -\frac{1}{78125}\e^8 +\frac{1}{390625}\e^9 -\frac{1}{9765675}\e^{10}\>.
\end{align}
%
We note the following. First, $z_2$ exactly solves the modified equation
\begin{align}
x^5-\e x-1\enskip +\frac{1}{25}\e^3 + \frac{3}{25}\e^4-\ldots + \frac{1}{9765625}\e^{10}=0 \label{starred}
\end{align}
which is $O(\e^3)$ different to the original. Second, the complete residual was computed rationally: there is no error in saying that $z_2=1+\sfrac{\e}{5}-\sfrac{\e^2}{25}$ solves equation \eqref{starred} exactly. Third, if $\e=1$ then $z_2=1+\sfrac{1}{5}-\sfrac{1}{25}=1.16$ exactly (or $1\sfrac{4}{25}$ if you prefer), and the residual is then $(\sfrac{29}{25})^5-\sfrac{29}{25}-1\doteq -0.059658$, showing that $1.16$ is the exact root of an equation about 6\% different to the original.

Something simple but importantly different to the usual treatment of perturbation methods has happened here. We have assessed the quality of the solution in an explicit fashion without concern for convergence issues or for the exact solution to $x^5-x-1=0$, which we term the reference problem. We use this term because its solution will be the reference solution. We can't call it the ``exact'' solution because $z_2$ is \textsl{also} an ``exact'' solution, namely to equation~\eqref{starred}.

Every numerical analyst and applied mathematician knows that this isn't the whole story---we need some evaluation or estimate of 
the effects of such perturbations of the problem. One effect is the difference between $z_2$ and $x$, the reference solution, and this is what people focus on. We believe this focus is sometimes excessive. The are other possible views. For instance, if the backward error is physically reasonable. 
%
As an example, if $\e=1$ and $z_2=1.16$ then $z_2$ exactly solves $y^5-y-a=0$ where $a\neq 1$ but rather $a\doteq 0.9403$. If the original equation was really $u^5-u-\alpha=0$ where $\alpha=1\pm 5\%$ we might be inclined to accept $z_2=1.16$ because, for all we know, we might have the true solution (even though we're outside the $\pm 5\%$ range, we're only just outside; and how confident are we in the $\pm5\%$, after all?). 

\subsubsection{Simple computer algebra solution}

The following Maple script can be used to solve this or similar problems $f(u;\e)=0$. Other computer algebra systems can also be used.
\begin{lstlisting}
# Perturbation solution of F(u;epsilon) = 0
restart; #top down execution should have a clean state to begin.
macro(e = epsilon); #saves typing "epsilon" every time.
F := z -> z^5-e*z-1;
# Zeroth order solution, by inspection, is
z := 1; #solve(eval(F(z),e=0),z);
A := coeff(series((D(F))(z), e, 1), e, 0); #A must not be 0 for regularity
N := 3; #number of terms
Delta := F(z); #initial residual, must be O(e)
# Now, the iteration:
for k to N do
	u := -coeff(series(Delta, e, k+1), e, k);
	z := z + u*e^k/A;
	Delta := F(z);
end do;
z;
series(Delta, e, N+3);
\end{lstlisting}
That code is a straightforward implementation of the general scheme presented in subsection \ref{genframe}. Its results, translated into \LaTeX\ and cleaned up a bit, are that
\begin{align} 
z = 1+\frac{1}{5}\e-\frac{1}{25}\e^2+\frac{1}{125}\e^3
\end{align}
and that the residual of this solution is
\begin{align}
\Delta = \frac{21}{3125}\e^5+O\left( \e^6 \right) \>.
\end{align}
With $N=3$, we get an extra order of accuracy as the next term in the series is zero, but this result is serendipitous.

\subsubsection{Systems of algebraic equations}\label{systems}
Regular perturbation for systems of equations using the framework from section \ref{genframe} is straightforward. We include an example to show some computer algebra and for completeness. Consider the following two equations in two unknowns:
\begin{align}
f_1(v_1,v_2) &=v_1^2+v_2^2 -1-\e v_1v_2 = 0\\
f_2(v_1,v_2) &= 25v_1v_2-12+2\e v_1 =0
\end{align}
When $\e=0$ these equations determine the intersections of a hyperbola with the unit circle. There are four such intersections: $(\sfrac{3}{5},\sfrac{4}{5}), (\sfrac{4}{5},\sfrac{3}{5}), (-\sfrac{3}{5},-\sfrac{4}{5})$ and $(-\sfrac{4}{5},-\sfrac{3}{5})$. The Jacobian matrix (which gives us the Fr\'echet derivative in the case of algebraic equations) is
\begin{align}
F_1(v) = \begin{bmatrix} \frac{\partial f_1}{\partial v_1} & \frac{\partial f_1}{\partial v_2} \\[.25cm] \frac{\partial f_2}{\partial v_1} & \frac{\partial f_2}{\partial v_2} \end{bmatrix} = \begin{bmatrix} 2v_1 & 2v_2 \\ 25 v_2 & 25 v_1\end{bmatrix} + O(\e)\>.
\end{align}
Taking for instance $u_0=[\sfrac{3}{5},\sfrac{4}{5}]^T$ we have
\begin{align}
A= F_1(u_0) = \begin{bmatrix} \sfrac{6}{5} & \sfrac{8}{5} \\ 20 & 15\end{bmatrix}\>.
\end{align}
Since $\det A=-14\neq 0$, $A$ is invertible and indeed
\begin{align}
A^{-1} = \begin{bmatrix} -\sfrac{15}{14} & \sfrac{4}{25} \\ \sfrac{10}{7} & -\sfrac{3}{35} \end{bmatrix}\>.
\end{align}
The residual of the zeroth order solution is
\begin{align}
\Delta_0 = F\left(\frac{3}{5},\frac{4}{5}\right) = \begin{bmatrix}-\sfrac{12}{25} \\ \sfrac{6}{5} \end{bmatrix}\>,
\end{align}
so $-[\e]\Delta_0 = [\sfrac{12}{25},-\sfrac{6}{5}]^T$. Therefore
\begin{align}
u_1 = \begin{bmatrix} u_{11} \\ u_{12}\end{bmatrix} = A^{-1}\begin{bmatrix}\sfrac{12}{25} \\ -\sfrac{6}{25}\end{bmatrix} = \begin{bmatrix} -\sfrac{114}{175} \\ \sfrac{138}{175}\end{bmatrix}
\end{align}
and $z_1=u_0+\e u_1$ is our improved solution:
\begin{align}
z_1 = \begin{bmatrix}\sfrac{3}{5} \\ \sfrac{4}{5} \end{bmatrix} + \e \begin{bmatrix} -\sfrac{114}{175} \\ \sfrac{138}{175}\end{bmatrix}\>.
\end{align}
To guard against slips, blunders, and bugs (some of those calculations were done by hand, and some were done in Sage on an Android phone) we compute
\begin{align}
\Delta_1 = F(z_1) = \e^2\begin{bmatrix}\sfrac{6702}{6125} \\ -\sfrac{17328}{1225}\end{bmatrix} + O\left(\e^3\right)\>.
\end{align}
That computation was done in Maple, completely independently. Initially it came out $O(\e)$ indicating that something was not right; tracking the error down we found a typo in the Maple data entry ($183$ was entered instead of $138$). Correcting that typo we find $\Delta_1=O(\e^2)$ as it should be. Here is the corrected Maple code:
\begin{lstlisting}
# Residual computation for a system of two equations
restart; #top down execution should have a clean state to begin.
macro(e = epsilon); #saves typing "epsilon" every time.
f1 := (v1,v2) -> v1^2 + v2^2 - 1 - e*v1*v2;
f2 := (v1,v2) -> 25*v1*v2 - 12 + 2*e*v1;
z11 := 3/5 + e*(-114/175);
z12 := 4/5 + e*138/175;
Delta11 := series( f1(z11,z12), e, 3);
Delta12 := series( f2(z11,z12), e, 3);
\end{lstlisting}
Just as for the scalar case, this process can be systematized and we give one way to do so in Maple, below. The code is not as pretty as the scalar case is, and one has to explicitly ``map'' the series function and the extraction of coefficients onto matrices and vectors, but this demonstrates feasibility.
\begin{lstlisting}
# Residual computation for a system of two equations
restart; #top down execution should have a clean state to begin.
macro(e = epsilon); #saves typing "epsilon" every time.
z := Vector(2,[3/5,4/5]); # z_0 = u_0
F := u -> Vector(2,
		[ u[1]^2 + u[2]^2 - 1 - e*u[1]*u[2],
		  25*u[1]*u[2] - 12 + 2*e*u[1] ]);
A := VectorCalculus[Jacobian](
		[ F([x,y])[1], F([x,y])[2]], [x,y] );
A := eval( A, [x=z[1], y=z[2], e=0] );
N := 3;
Delta := F(z);
for k to N do
	u := map(t -> -coeff( t, e, k ),
		 map( series, Delta, e, k+1 )
		);
	z := z + LinearAlgebra[LinearSolve]( A, u )*e^k;
	Delta := F( z );
end do:
z;
map( series, Delta, e , N+2 );
\end{lstlisting}
This code computes $z_3$ correctly and gives a residual of $O(\e^4)$. From the backward error point of view, this code finds the intersection of curves that differ from the specified ones by terms of $O(\e^4)$. In the next section, we show  a way to use a built-in feature of Maple to do the same thing with less human labour.

\subsubsection{The Davidenko equation}
Maple has a built-in facility for solving differential equations in series that (at the time of writing) is superior to its built-in facility for solving algebraic equations in series, because the latter can only handle scalar equations. This may change in the future, but it may not because there is the following simple workaround. To solve
\begin{align}
F(u;\e)=0
\end{align}
for a function $u(\e)$ expressed as a series, simply differentiate to get
\begin{align}
D_1(F)(u,\e)\frac{du}{d\e} + D_2(F)(u,\e)=0\>.
\end{align}
Boyd \cite{Boyd(2014)} calls this the Davidenko equation. If we solve this in Taylor series with the initial condition $u(0)=u_0$, we have our perturbation series. Notice that what we were calling $A=[\e^0]F_1(u_0)$ occurs here as $D_1(F)(u_0,0)$ and this needs to be nonsingular to be solved as an ordinary differential equation; if $\mathrm{rank}(D_1(F)(u_0,0))<n$ then this is in fact a nontrivial differential algebraic equation that Maple may still be able to solve using advanced techniques (see, e.g., \cite{Avrachenkov(2013)}). Let us just show a simple case here:
\begin{lstlisting}
# Residual computation for a system of two equations
restart; #top down execution should have a clean state to begin.
macro(e = epsilon); #saves typing "epsilon" every time.
Order := 4;
z := Vector(2,[3/5,4/5]); # z_0 = u_0
F := u -> Vector(2,
		[ u[1]^2 + u[2]^2 - 1 - e*u[1]*u[2],
		  25*u[1]*u[2] - 12 + 2*e*u[1] ]);
Zer := F( [ x(e), y(e) ] ); #This asks for F evaluated at functions x(e)
# and y(e) that are yet unspecified.
diffeqs := { diff( Zer[1], e ), diff( Zer[2], e ) }; #This creates a set
	# of two differential equations, one from each component of F.
	# Each equation will contain both dx/de and dy/de.
iniconds := { x(0) = z[1] , y(0) = z(2) };
sol := dsolve( diffeqs union iniconds, {x(e), y(e)}, type=series );
Delta := eval( F( [x(e), y(e)] ), map(convert, sol, polynom) ):
map( series, Delta, e, Order+2 );
\end{lstlisting}
This generates (to the specified value of the order, namely, \verb|Order=4|) the solution
\begin{align}
 x(\e) &=\frac{3}{5}-\frac{114}{175}\e+\frac{119577}{42875}\e^2-\frac{43543632}{2100875}\e^3\\
y(\e) &=\frac{4}{5}+\frac{138}{175}\e-\frac{119004}{42875}\e^2+\frac{43245168}{2100875}\e^3\>,
\end{align}
whose residual is $O(\e^4)$. Internally, Maple uses its own algorithms, which occasionally get improved as algorithmic knowledge advances.

\subsection{Puiseux series}\label{Puiseux}
Puiseux series 
 are simply Taylor series or Laurent series with fractional powers. A standard example is
\begin{align}
\sin\sqrt{x} = x^{\sfrac{1}{2}} - \frac{1}{3!}x^{\sfrac{3}{2}} + \frac{1}{5!}x^{\sfrac{5}{2}}+\cdots
\end{align}
A simple change of variable (e.g. $t=\sqrt{x}$ so $x=t^2$) is enough to convert to Taylor series. Once the appropriate power $n$ is known for $\e=\mu^n$, perturbation by Puiseux expansion reduces to computations similar to those we've seen already.
For instance, had we chosen to embed $u^5-u-1$ in the family $u^5-\e(u+1)$ (which is somehow conjugate to the family of the last section), then because the equation becomes $u^5=0$ when $\e=0$ we see that we have a five-fold root to perturb, and we thus suspect we will need Puiseux series.

For scalar equations, there are built-in facilities in Maple for Puiseux series, which gives yet another way in Maple to solve scalar algebraic equations perturbatively. One can use the \texttt{RootOf} construct to do so as follows:
\begin{lstlisting}
restart;
macro(e = epsilon);
Order := 2;
alias(alpha = RootOf(z^5-1, z));
f := u -> u^5 - e*(u+1);
z := convert( series( RootOf(f(u),u), e), polynom);
Delta := series(f(z), e, Order+2):
map(simplify, Delta);
\end{lstlisting}
This yields
\begin{align}
z = \alpha\e^{\sfrac{1}{5}}+\frac{1}{5}\alpha^2\e^{\sfrac{2}{5}} -\frac{1}{25}\alpha^3\e^{\sfrac{3}{5}} +\frac{1}{125}\alpha^4\e^{\sfrac{4}{5}} - \frac{21}{15626}\alpha \e^{\sfrac{6}{5}} \>.
\end{align}
This series describes all paths, accurately for small $\e$. Note that the command
\begin{lstlisting}
alias(alpha = RootOf(u^5-1,u))
\end{lstlisting}
 is a way to tell Maple that $\alpha$ represents a fixed fifth root of unity. Exactly which fixed root can be deferred till later. Working instead with the default value for the environment variable \texttt{Order}, namely \texttt{Order := 6}, gets us a longer series for $z$ containing terms up to $\e^{\sfrac{29}{5}}$ but not $\e^{\sfrac{30}{5}}=\e^6$. Putting the resulting $z_6$ back into $f(u)$ we get a residual
\begin{align}
\Delta_6 = f(z_6) = \frac{23927804441356816}{14551915228366851806640625}\e^7 + O(\e^8)
\end{align}
Thus we expect that for small $\e$ the residual will be quite small. For instance, with $\e=1$ the exact residual is, for $\alpha=1$, $\Delta_6=1.2\cdot 10^{-9}$. This tells us that this approximation ought to get us quite accurate roots, and indeed we do.

We conclude this discussion with two remarks. The first is that by a discriminant analysis as we describe in section \ref{SingPert}, we find that the nearest singularity is at $\e=\sfrac{3125}{256}$, and so we expect this series to actually converge for $\e=1$. Again, this fact was not used in our analysis above. Secondly, we could have used the \verb|series/RootOf| technique to do both the regular perturbation in subsection \ref{RegularPert} or the singular one we will do in subsection \ref{SingPert}. The Maple commands are quite similar:
\begin{lstlisting}
series(RootOf(u^5-e*u-1,u),e);
\end{lstlisting}
and
\begin{lstlisting}
series(RootOf(e*u^5-u-1,u),e);
\end{lstlisting}
However, in both cases only the real root is expanded. Some ``Maple art'' (that one of us more readily characterizes as black magic) can be used to complete the computation, but the previous code (both the loop and the Davidenko equation) are easier to generalize. Making the \texttt{dsolve/series} code for the Davidenko equation work in the case of Puiseux series requires a preliminary scaling.

\subsection{Singular perturbation}\label{SingPert}
%
%
%
%
%

Suppose that instead of embedding $u^5-u-1=0$ in the regular family we used in the previous section, we had used $\e u^5-u-1=0$. If we run our previous Maple programs, we find that the zeroth order solution is unique, and $z_0=-1$. The Fr\'echet derivative is $-1$ to $O(\e)$, and so $u_{n+1} = [\e^{n+1}]\Delta_n$ for all $n\geq 0$. We find, for instance,
\begin{align}
z_7 = -1-\e -5\e^2 - 35\e^3 -285\e^4 -2530\e^5 -23751\e^6 -231880\e^7
\end{align}
which has residual $\Delta_7 = O(\e^8)$ but with a larger integer as the constant hidden in that $O$ symbol. For $\e=0.2$, the value of $z_7$ becomes \begin{align}z_7\doteq -7.4337280\end{align} while $\Delta_7=-4533.64404$, which is not small at all. Thus we have no evidence this perturbation solution is any good: we have the exact solution to $u^5-0.2 u-1=-4533.64404$ or $u^5-0.2 u+4532.64404=0$, probably not what was intended (and if it was, it would be a colossal fluke). Note that we do not need to know a reference value of a root of $u^5-0.2u-1$ to determine this. 
Trying a smaller $\e$, we find that if $\e=0.05$ we have $z_7\doteq -1.07$ and $\Delta_7\doteq -1.2\cdot 10^{-4}$. This means $z_7$ is an exact root of $u^5-0.05 u-1.00012$; which may very well be what we want.

The following remark is not really germane to the method but it's interesting. Taking the discriminant with respect to $u$, i.e., the resultant of $f$ and $\sfrac{\partial f}{\partial u}$, we find $\mathrm{discrim}(f) = \e^3(3125\e-256)$. Thus $f$ will have multiple roots if $\e=0$ (there are 4 multiple roots at infinity) or if $\e = \sfrac{256}{3125}=0.08192$. Thus our perturbation expansion can be expected to diverge\footnote{A separate analysis leads  to the identification of $u_k = \frac{1}{5k+1}\binom{5k+1}{k}$ (via \cite{OEIS}). The ratio test confirms that the series converges for $|\e|<\sfrac{256}{3125}$, and diverges if $\e = \sfrac{256}{3125}$.} for $\e\geq 0.08192$. What happens to $z_7$ if $\e=\sfrac{256}{3125}$? $z_7\doteq -1.1698$ and $\Delta_7=-9.65\cdot10^{-3}$, so we have an exact solution for $u^5-\sfrac{256}{3125}u-1.00965$; this is not bad. The reference double root is $-1.25$, about $0.1$ away, although this fact was not used in the previous discussion. 

But this computation, valid as it is, only found one root out of five, and then only for sufficiently small $\e$. We now turn to the roots that go to infinity as $\e\to0$. Preliminary investigation from similar to that of subsection \ref{Puiseux} shows that it is convenient to replace $\e$ by $\mu^4$. 
Many singular perturbation problems including this one can be turned into regular ones by rescaling. Putting $u=\sfrac{y}{\mu}$, we get
\begin{align}
\mu^4\left(\frac{y}{\mu}\right)^5-\frac{y}{\mu}-1=0\>,
\end{align}
which reduces to
\begin{align}
y^5-y-\mu=0\>.
\end{align}
This is now regular in $\mu$. The zeroth order the equation is $y(y^4-1)=0$ and the root $y=0$ just recovers the regular series previously attained; so we let $\alpha$ be a root of $y^4-1$, i.e., $\alpha\in\{1,-1,i,-i\}$. A very similar Maple program (to either of the previous two) gives
\begin{align}
y_5= \alpha +\frac{1}{4}\mu - \frac{5}{32}\alpha^3\mu^2 +\frac{5}{32}\alpha^2\mu^3-\frac{385}{2048}\alpha\mu^4 + \frac{1}{4}\mu^5
\end{align}
so our appoximate solution is $\sfrac{y_5}{\mu}$ or
\begin{align}
z_5 = \frac{\alpha}{\mu}+\frac{1}{4}-\frac{5}{32}\alpha^3\mu^2-\frac{385}{2048}\alpha\mu^3+\frac{1}{4}\mu^4
\end{align}
which has residual \textsl{in the original equation}
\begin{align}
\Delta_5 = \mu^4 z^5 -z-1= \frac{23205}{16384}\alpha^3\mu^5 - \frac{21255}{65536}\alpha^2\mu^6 +O(\mu^7)\>.\label{residorigeq}
\end{align}
That is, $z_5$ exactly solves $\mu^4u^5-u-1-\sfrac{23205}{16384}\>\alpha^2\mu^5=O(\mu^6)$ instead of the one we had wanted to solve. This differs from the original by $O(|\e|^{\sfrac{5}{4}})$, and for small enough $\e$ this may suffice.

\paragraph{Optimal backward error}
Interestingly enough, we can do better. The residual is only one kind of backward error. Taking the lead from the Oettli-Prager theorem \cite[chap.~6]{CorlessFillion(2013)}, we look for equations of the form
\begin{align}
\left(\mu^4 +\sum_{j=10}^{15} a_j\mu^j\right) u^5 - u -1
\end{align}
for which $z_5$ is a better solution yet. Simply equating coefficients of the residual
\begin{align}
\tilde{\Delta}_5 = \left(\mu^4+\sum_{j=10}^{15}a_j\mu^j\right)z_5^5-z_5-1
\end{align}
to zero, we find
\begin{align}
(\mu^4 - \frac{23205}{16384}\alpha^2\mu^{10}+ \frac{2145}{1024}\alpha\mu^{11})z_5^5 -z_5-1 = \frac{12165535425}{1073741824}\alpha\mu^{11}+O(\mu^{12})
\end{align}
and thus $z_5$ solves an equation that is $O(\mu^{\sfrac{10}{4}})=O(\e^{\sfrac{5}{2}})$ close to the original, not just an equation \eqref{residorigeq} that is $O(\mu^6)=O(|\e|^{\sfrac{5}{4}})$. This is a superior explanation of the quality of $z_5$. 
This was obtained with the following Maple code:
\begin{lstlisting}
# Perturbation solution of F(u;epsilon) = 0
restart;
e := mu^4;
Forig := z -> e*z^5 - z - 1;
F := y -> y^5 - y - mu;
# Zeroth order solution, by inspection:
alias(alpha = RootOf(Z^4-1, Z));
y := alpha;
A := coeff( series( (D(F))(y), mu, 1), mu, 0);
A := simplify(A);
N := 5;
Delta := simplify( F(y) );
for k to N do
	u := -coeff( series(Delta, mu, k+1), mu, k);
	y := y+u*mu^k/A;
	Delta := simplify( F(y) );
end do:
y;
series(Delta, mu, N+3);
M := 5+2*N;
modified := u -> (mu^4 + add(a[j]*mu^j, j = 5+N..M))*u^5 - u - 1;
z := series(y/mu, mu, N+1);
zer := series(modified(z), mu, M+1);
eqs := [seq(simplify(coeff(zer, mu, k)), k = N .. M-5)];
sol := solve(eqs, [seq(a[j], j = 5+N .. M)]);
perteq := eval(modified(U), sol[1]);
newresid := eval(perteq, U = z):
map(simplify, series(newresid, mu, M+2));
\end{lstlisting}
%
%
Computing to higher orders (see the worksheet) gives e.g. that $z_8$ is the exact solution to an equation that differs by $O(\mu^{13})$ from the original, or better than $O(\e^3)$. This in spite of the fact that the basic residual $\Delta_8=O(\e^{9/4})$, only slightly better than $O(\e^2)$.

We will see other examples of improved backward error over residual for singularly-perturbed problems. In retrospect it's not so surprising, or shouldn't have been: singular problems are sensitive to changes in the leading term, and so it takes less effort to match a given solution.

\subsection{Perturbing all roots at once}
The preceding analysis found a nearby equation for each root independently; this might suffice, but there are circumstances in which it might not. Perhaps we want a ``nearby'' equation satisfied by all roots at once. Sadly this is more difficult, and in general may not be possible. But it is possible for the example we've considered and we demonstrate how the backward error is used in such a case. Let
\begin{align}
\zeta_1 &= z_5(1)= \frac{1}{\mu}+\frac{1}{4} -\frac{5}{32}\mu-\frac{385}{2048}\mu^3+\frac{1}{4}\mu^4\\
\zeta_2 &= z_5(-1) = -\frac{1}{\mu}+\frac{1}{4}-\frac{5}{32}\mu +\frac{385}{2048}\mu^3 + \frac{1}{4}\mu^4\\
\zeta_3 &= z_5(i) = \frac{i}{\mu}+\frac{1}{4} + \frac{5}{32}\mu -\frac{385i}{2048}\mu^3 + \frac{1}{4}\mu^4\\
\zeta_4 &= z_5(-i) = -\frac{i}{\mu}+\frac{1}{4} + \frac{5}{32}\mu + \frac{385}{2048}\mu^3 + \frac{1}{4}\mu^4\\
\zeta_5 &= z_5 = -1-\mu^4-5\mu^8 ,
\end{align}
$\zeta_5$ is the regular root we have found first in the previous subsection. Now put 
\begin{align}
\tilde{p}(x)=\mu^4(x-\zeta_1)(x-\zeta_2)(x-\zeta_3)(x-\zeta_4)(x-\zeta_5)
\end{align}
and expand it. The result, by Maple, is
\begin{multline}
\mu^4x^5-5\mu^{12}x^4 + \left(\frac{23205}{16384}\mu^8+\frac{45}{8}\mu^{12}\right)x^3
 -\left(\frac{5435}{32768}\mu^8+\frac{195697915}{33554432}\mu^{12}\right)x^2  \\
 + \left( \frac{2575665}{2097152}\mu^8+\frac{5696429035}{1073741824}\mu^{12}-1 \right)x + 
\frac{8453745}{2097152}\mu^8  -\frac{5355037365}{1073741824}\mu^{12}-1
\end{multline}
which equals
\begin{align}
\e x^5 -x-1-5\e^3 x^4 + (\frac{23205}{16384}\e^2+\frac{45}{8}\e^3)x^3 - (\frac{5435}{32768}\e^2+\cdots)x^2+O(\e^2)
\end{align}
As we see, this equation is remarkably close to the original, although we see changes in all the coefficients. The backward error is $O(\mu^8)$, i.e., $O(\e^2)$. Thus for algebraic equations it's possible to talk about simultaneous backward error.

\subsection{A hyperasymptotic example}
In \cite[sect.~15.3, pp.~285-288]{Boyd(2014)}, Boyd takes up the perturbation series expansion of the root near $-1$ of
\begin{align}
f(x,\e)=1+x+\e \mathrm{sech}\left(\frac{x}{\e}\right) = 0\>,
\end{align}
a problem he took from \cite[p.~22]{Holmes(1995)}. After computing the desired expansion using a two-variable technique, Boyd then sketches an alternative approach suggested by one of us (based on \cite{Corless(1996)}), namely to use the Lambert $W$ function. Unfortunately, there are a number of sign errors in Boyd's equation (15.28). We take the opportunity here to offer a correction, together with a residual-based analysis that confirms the validity of the correction. First, the erroneous formula: Boyd has
\begin{align}
z_0 = \frac{W(-2e^{\sfrac{1}{\e}})\e-1}{\e}
\end{align}
and $x_0=-\e z_0$, so allegedly $x_0=1-\e W(-2\e^{\sfrac{1}{\e}})$. This can't be right: as $\e\to0^+$, $e^{\sfrac{1}{\e}}\to\infty$ and the argument to $W$ is negative and large; but $W$ is real only if its argument is between $-e^{-1}$ and $0$, if it's negative at all. We claim that the correct formula is
\begin{align}
x_0 = -1-\e W(2e^{-\sfrac{1}{\e}}) \label{star}
\end{align}
which shows that the errors in Boyd's equation (15.28) are explainable as trivial. Indeed, Boyd's derivation is correct up to the last step; rather than fill in the algebraic details of the derivation of formula~\eqref{star}, we here verify that it works by computing the residual:
\begin{align}
\Delta_0 = 1+x_0 + \e \mathrm{sech}\left(\frac{x_0}{\e}\right).
\end{align}
For notational simplicity, we will omit the argument to the Lambert $W$ function and just write $W$ for $W(2e^{-\sfrac{1}{\e}})$. Then, note that $\mathrm{sech}(\sfrac{x_0}{\e}) = \mathrm{sech}(\sfrac{1+\e W}{\e})$ since each $\mathrm{sech}$ is even, and that
\begin{align}
\mathrm{sech}\left(\frac{x_0}{\e}\right) = \frac{2}{\displaystyle e^{\sfrac{x_0}{\e}}+e^{-\sfrac{x_0}{\e}}} = \frac{1}{\displaystyle e^{(\sfrac{1}{\e}) +W}+e^{-\sfrac{1}{\e}-W}}\>.
\end{align}
Now, by definition,
\begin{align}
We^W = 2e^{-\sfrac{1}{\e}}
\end{align}
and thus we obtain
\begin{align}
e^W = \frac{2e^{-\sfrac{1}{\e}}}{W} \qquad \textrm{and} \qquad e^{-W} = \frac{We^{\sfrac{1}{\e}}}{2}\>.
\end{align}
It follows that
\begin{align}
\mathrm{sech}\left(\frac{x_0}{\e}\right) = \frac{2}{\displaystyle \sfrac{2}{W}+\sfrac{W}{2}} = \frac{W}{\displaystyle 1+\sfrac{W^2}{4}}\>, \label{sechW}
\end{align}
and hence the residual is
\begin{align}
\Delta_0 &= 1+(-1-\e W)+\e \frac{W}{\displaystyle 1+\sfrac{W^2}{4}} 
= \frac{\displaystyle -\e W(1+\sfrac{W^2}{4}) + \e W}{\displaystyle 1+\sfrac{W^2}{4} } \\
&= \frac{\displaystyle -\sfrac{\e W^3}{4}}{\displaystyle 1+\sfrac{W^2}{4}} = \frac{-\e W^3}{4+ W^2} \nonumber \>.
\end{align}
Now $W= W(2e^{-1/\e})$ and as $\e\to 0^+$, $2e^{-1/\e}\to 0$ rapidly; since the Taylor series for $W(z)$ starts as $W(z)= z-z^2+\frac{3}{2}z^3+\ldots$, we have that $W(2e^{-\sfrac{1}{\e}})\sim 2e^{-\sfrac{1}{\e}}$ and therefore
\begin{align}
\Delta_0 = -\e 2e^{-\sfrac{3}{\e}}+O(e^{-\sfrac{5}{\e}})\>.
\end{align}
We see that this residual is very small indeed. But we can say even more. Boyd leaves us the exercise of computing higher order terms; here is our solution to the exercise. A Newton correction would give us
\begin{align}
x_1 = x_0 - \frac{f(x_0)}{f'(x_0)}
\end{align}
and we have already computed $f(x_0)=\Delta_0$. What is $f'(x_0)$? Since $f(x) = 1+x+\e\mathrm{sech}(\sfrac{x}{\e})$, this derivative is
\begin{align}
f'(x) = 1-\mathrm{sech}\left(\frac{x}{\e}\right)\mathrm{tanh}\left(\frac{x}{\e}\right)\>.
\end{align}
Simplifying similarly to equation \eqref{sechW}, we obtain
\begin{align}
\mathrm{tanh}\left(\frac{x_0}{\e}\right) = \frac{e^{1/\e +W} - e^{-1/\e-W}}{e^{1/\e+W}+e^{-1/\e+W}} = \frac{\frac{2}{W}-\frac{W}{2}}{\frac{2}{W}+\frac{W}{2}} = \frac{4-W^2}{4+W^2}\>.
\end{align}
Thus
\begin{align}
f'(x_0) &= 1-\mathrm{sech}\left(\frac{x_0}{\e}\right)\mathrm{tanh}\left(\frac{x_0}{\e}\right)
= 1- \frac{\displaystyle W(1-\sfrac{W^2}{4})}{\displaystyle (1+\sfrac{W^2}{4})^2}\>.
\end{align}
It follows that
\begin{align}
x_1 &= x_0 - \frac{\Delta_0}{f'(x_0)} 
= -1-\e W+\frac{\displaystyle\sfrac{\e W^3}{4+W^2}}{\displaystyle 1- \frac{W(1-\sfrac{W^2}{4})}{(1+\sfrac{W^2}{4})^2}} \\
&= -1 -\e W+ \frac{\e W^3(4+W^2)}{16-16W+8W^2+4W^3+W^4}\\
&= -1-\e W+ \frac{\e}{4}W^3+\frac{\e}{4}W^4+\frac{3}{16}\e W^5-\frac{11}{64}\e W^6+O(W^7)
\end{align}
Finally, the residual of $x_1$ is
\begin{align}
\Delta_1 = 4\e e^{\sfrac{7}{\e}}+O(\e e^{-\sfrac{8}{\e}})\>. \label{Newtresid}
\end{align}
We thus see an example of the use of $f'(x_0)$ instead of just $A$, as discussed in section \ref{genframe}, to approximately double the number of correct terms in the approximation.

This analysis can be implemented in Maple as follows:
\begin{lstlisting}
restart;
with(MultiSeries);
macro(e = epsilon);
alias(W = LambertW);
f := x -> 1 + x + e*sech(x/e);
df := D(f);
x[0] := -1 - e*W(2*exp(-1/e));
Delta[0] := f( x[0] );
series(Delta[0], e, 3);
x[1] := x[0] - Delta[0]/df(x[0]);
Delta[1] := f(x[1]);
s := multiseries(x[1],e=0);
scale := SeriesInfo[Scale](s);
multiseries(x[1],scale,3);
multiseries(Delta[1],scale,5);
# In what follows we have substituted expressions in W for sech and tanh
# since Maple couldn't simplify the expression well.
restart;
macro(e = epsilon);
x[1] := -1-e*W+e*W^3/((4+W^2)*(1-W*(1-(1/4)*W^2)/(1+(1/4)*W^2)^2));
change := factor(x[1]+1+e*W);
series(change,W=0,8);
\end{lstlisting}
Note that we had to use the MultiSeries package \cite{Salvy(2010)} to expand the series in equation \eqref{Newtresid}, for understanding how accurate $z_2$ was. $z_2$ is slightly more lacunary than the two-variable expansion in \cite{Boyd(2014)}, because we have a zero coefficient for $W^2$.

\section{Divergent Asymptotic Series}
Before we begin, a note about the section title: some authors give the impression that the word ``asymptotic'' is used \textsl{only} for divergent series, 
 and so the title might seem redundant. But the proper definition of an asymptotic series can include convergent series (see, e.g., \cite{Bruijn(1981)}), as it means that the relevant limit is not as the number of terms $N$ goes to infinity, but rather as the variable in question (be it \e, or $x$, or whatever) approaches a distinguished point (be it 0, or infinity, or whatever). In this sense, an asymptotic series might diverge as $N$ goes to infinity, or it might converge, but typically we don't care. We concentrate in this section on divergent asymptotic series.

 Beginning students are often confused when they learn the usual ``rule of thumb'' for optimal accuracy when using divergent asymptotic series, namely to truncate the series \textsl{before} adding in the smallest (magnitude) term. This rule is usually motivated by an analogy with \textsl{convergent} alternating series, where the error is less than the magnitude of the first term neglected. But why should this work (if it does) for divergent series?

The answer we present in this section isn't as clear-cut as we would like, but nonetheless we find it explanatory. Perhaps you and your students will, too. The basis for the answer is that one can measure the residual $\Delta$ that arises on truncating the series at, say, $M$ terms, and choose $M$ to minimize the residual. Since the forward error is bounded by the condition number times the size of the residual, by minimizing $\|\Delta\|$ one minimizes a bound on the forward error. It often turns out that this method gives the same $M$ as the rule of thumb, though not always.

An example may clarify this. We use the large-$x$ asymptotics of $J_0(x)$, the zeroth-order Bessel function of the first kind. In \cite[section 10.17(i)]{NIST:DLMF}, we find the following asymptotic series, which is attributed to Hankel:
\begin{align}
J_0(x) = \left(\frac{2}{\pi x}\right)^{\sfrac{1}{2}}\left( A(x)\cos\left(x-\frac{\pi}{4}\right)-B(x)\sin\left(x-\frac{\pi}{4}\right)\right)
\end{align}
where
\begin{align}
A(x) = \sum_{k\geq 0} \frac{a_{2k}}{x^{2k}} \qquad \textrm{and}\qquad 
B(x) = \sum_{k\geq 0} \frac{a_{2k+1}}{x^{2k+1}} \label{twoseries}
\end{align}
and where
\begin{align}
a_0 &= 1 \nonumber\\
a_k &= \frac{(-1)^k }{k! 8^k}\prod_{j=1}^k (2j-1)^2\>.
\end{align}
For the first few $a_k$s, we get
\begin{align}
a_0=1, a_1 = -\frac{1}{8}, a_2= -\frac{9}{128}, a_3 = \frac{75}{1024}\>,
\end{align}
and so on. The ratio test immediately shows the two series \eqref{twoseries} diverge for all finite~$x$.

Luckily, we always have to truncate anyway, and if we do, the forward errors get arbitrarily small so long as we take $x$ arbitrarily large. Because the Bessel functions are so well-studied, we have alternative methods for computation, for instance
\begin{align}
J_0(x) = \frac{1}{\pi}\int_0^\pi \cos(x\sin\theta)d\theta
\end{align}
which, given $x$, can be evaluated numerically (although it's ill-conditioned in a relative sense near any zero of $J_0(x)$). So we can directly compute the forward error. 
But let's pretend that we can't. We have the asymptotic series, and not much more. Or course we have to have a defining equation---Bessel's differential equation
\begin{align}
x^2y''+xy'+x^2y=0
\end{align}
with the appropriate normalizations at $\infty$. We look at
\begin{align}
y_{N,M} = \left(\frac{2}{\pi x}\right)^{\sfrac{1}{2}}A_N(x)\cos\left(x-\frac{\pi}{4}\right)-\frac{2}{\pi x}B_M(x)\cos\left(x-\frac{\pi}{4}\right)
\end{align}
where
\begin{align}
A_N(x) = \sum_{k=0}^N \frac{a_{2k}}{x^{2k}}\qquad \textrm{and}\qquad 
B_M(x) = \sum_{k=0}^M \frac{a_{2k+1}}{x^{2k+1}}\>.
\end{align}
Inspection shows that there are only two cases that matter: when we end on an even term $a_{2k}$ or on an odd term $a_{2k+1}$. The first terms omitted will be odd and even. A little work shows that the residual
\begin{align}
\Delta = x^2y''_{N,M} + x y'_{N,M} + x^2y_{N,M}
\end{align}
is just
\begin{align}
\frac{\displaystyle (k+\sfrac{1}{2})^2 a_k}{\displaystyle x^{k+\sfrac{1}{2}}} \cdot \left\{ \begin{array}{c} \cos(x-\sfrac{\pi}{4})\\ \sin(x-\sfrac{\pi}{4})\end{array}\right\}
\end{align}
if the final term \textsl{kept}, odd or even, is $a_k$. If even, then multiply by $\cos(x-\sfrac{\pi}{4})$; if odd, then $\sin(x-\sfrac{\pi}{4})$.

Let's pause a moment. The algebra to show this is a bit finicky but not hard (the equation is, after all, linear). This end result s an extremely simple (and exact!) formula for $\Delta$. The finite series $y_{N,M}$ is then the exact solution to
\begin{align}
x^2y''+xy'+xy &= \Delta\\
&=  \frac{\displaystyle (k+\sfrac{1}{2})^2 a_k}{x^{k+\sfrac{1}{2}}} \cdot \left\{ \begin{array}{c} \cos(x-\frac{\pi}{4})\\ \sin(x-\frac{\pi}{4})\end{array}\right\}
\end{align}
and, provided $x$ is large enough, this is only a small perturbation of Bessel's equation. In many modelling situations, such a small perturbation may be of direct physical significance, and we'd be done. Here, though, Bessel's equation typically arises as an intermediate step, after separation of variables, say. Hence one might be interested in the forward error. By the theory of Green's functions, we may express this as
\begin{align}
J_0(x) - y_{N,M}(x) = \int_x^\infty K(x,\xi)\Delta(\xi)d\xi
\end{align}
for a suitable kernel $K(x,\xi)$. The obvious conclusion is that if $\Delta$ is small then so will $J_0(x)-y_{N,M}(x)$; but $K(x,\xi)$ will have some effect, possibly amplifying the effects of $\Delta$, or perhaps even damping its effects. Hence, the connection is indirect.

To have an error in $\Delta$ of at most $\e$, we must have
\begin{align}
\left(k+\frac{1}{2}\right)^2\frac{|a_k|}{x^{k+\sfrac{1}{2}}}\leq\e
\end{align}
(remember, $x>0$). This will happen only if
\begin{align}
x\geq \left(\left(k+\frac{1}{2}\right)^2 \frac{|a_k|}{\e}\right)^{2/(2k+1)}
\end{align}
and this, for fixed $k$, goes to $\infty$ as $\e\to0$. 
Alternatively, we may ask which $k$, for a fixed $x$, minimizes
\begin{align}
\left(k+\frac{1}{2}\right)^2\frac{|a_k|}{x^{k+\sfrac{1}{2}}}
\end{align}
and this answers the truncation question in a rational way. In this particular case, minimizing $\|\Delta\|$ doesn't necessarily minimize the forward error (although, it's close). For $x=2.3$, for instance, the sequence $(k+\sfrac{1}{2})^2|a_k|x^{-k-\sfrac{1}{2}}$ is (no $\sqrt{\sfrac{2}{\pi}}$)
\begin{align}
\begin{array}{ccccccc}
k & 0 & 1 & 2 & 3 & 4 & 5\\
A_k & 0.165 & 0.081 & 0.055 & 0.049 & 0.054 & 0.070
\end{array}
\end{align}
The clear winner seems to be $k=3$. This suggests that for $x=2.3$, the best series to take is
\begin{align}
y_3 = \left(\frac{2}{\pi x}\right)^{\sfrac{1}{2}} \left( \left(1-\frac{9}{128x^2}\right)\cos\left(x-\frac{\pi}{4}\right) + \left(\frac{1}{8x}-\frac{75}{1024x^3}\right)\sin\left(x-\frac{\pi}{4}\right)\right)\>.
\end{align}
This gives $5.454\cdot 10^{-2}$ for $x=2.3$. But the cosine versus sine plays a role, here: $\cos(2.3-\sfrac{\pi}{4})\doteq 0.056$ while $\sin(2.3-\sfrac{\pi}{4})\doteq0.998$, so we should have included this. When we do, the estimates for $\Delta_0,\Delta_2$ and $\Delta_4$ are all significantly reduced---and this changes our selection, and makes $k=4$ the right choice; $\Delta_6>\Delta_4$ as well (either way). But the influence of the integral is mollifying. 
Comparing to a better answer (computers via the integral formula) $0.0555398$, we see that the error is about $8.8\cdot 10^{-4}$ whereas $((4+\sfrac{1}{2})^2a_4/2.3^{4+\sfrac{1}{2}})\cos(2.3-\sfrac{\pi}{4})$ is $3.06\cdot 10^{-3}$; hence the residual overestimates the error slightly.

How does the rule of thumb do? The first term that is neglected here is $(\sfrac{1}{x})^{\sfrac{1}{2}}a_5x^{-5}\sin(x-\sfrac{\pi}{4})$ which is $\sim2.3\cdot 10^{-3}$ apart from the $(\sfrac{2}{\pi})^{\sfrac{1}{2}}=0.797$ factor, so about $1.86\cdot10^{-3}$. The \textsl{next} term is, however, $(\sfrac{2}{\pi x})^{\sfrac{1}{2}}a_6x^{-6}\cos(x-\sfrac{\pi}{4})\doteq -1.14\cdot 10^{-4}$ which is smaller yet, suggesting that we should keep the $a_5$ term. 
But we shouldn't. Stopping with $a_4$ gives a better answer, just as the residual suggests that it should.

We emphasize that this is only a slightly more rational rule of thumb, because minimizing $\|\Delta\|$ only minimizes a bound on the forward error, not the forward error itself. Still, we have not seen this discussed in the literature before. A final comment is that the defining equation and its scale, define also the scale for what's a ``small'' residual.

So, a justification for the ``rule of thumb'' would be as follows. In our general scheme,
\begin{align}
Au_{n+1} = -[\e^{n+1}]\Delta_n
\end{align}
and thus, loosely speaking,
\begin{align}
u_{n+1} \sim -A^{-1}\Delta_n + O(\e^{n+1})\>.
\end{align}
Thus, if we stop when $u_{n+1}$ is smallest, this would tend to happen at the same integer $n$ that $\Delta_n$ was smallest.

This isn't going to be always true. For instance, if $A$ is a matrix with largest singular value $\sigma_1$ and smallest $\sigma_N>0$, with associated vectors $\hat{u}_k$ and $\hat{v}_k$, so that
\begin{align}
A\hat{v}_k = \sigma_k\hat{u}_k\>.
\end{align}
Then, if $u_{n+1}$ is like $\hat{v}_1$ then $\Delta_n$ will be like $\sigma\hat{u}_1$, which can be substantially larger; contrariwise, if $u_{n+1}$ is like $\hat{v}_N$ then $A\hat{v}_N=\sigma_N\hat{u}_N$ and $\Delta_n$ can be substantially smaller. The point is that directions of $\Delta_n$ can change between steps in the perturbation expansion; we thus expect correlation but not identity.

\section{Initial-Value problems}

BEA has successfully been applied to the \textsl{numerical} solution of differential equations for a long time, now. Examples include the works of Enright since the 1980s, e.g., \cite{Enright(1989)b,Enright(1989)a}, and indeed the Lanczos $\tau$-method is yet older~\cite{Lanczos(1988)}. It was pointed out in \cite{Corless(1992)} and \cite{Corless(1993)b} that BEA could be used for perturbation and other series solutions of differential equations, also. We here display several examples illustrating this fact. We use regular expansion, matched asymptotic expansions, the renormalization group method, and the method of multiple scales.

\subsection{Duffing's Equation}
This proposed way of interpreting solutions obtained by perturbation methods has interesting advantages for the analysis of series solutions to differential equations. Consider for example an unforced weakly nonlinear Duffing oscillator, which we take from \cite{Bender(1978)}:
\begin{align}
y''+y+\varepsilon y^3=0 \label{Duffing}
\end{align}
with initial conditions $y(0)=1$ and $y'(0)=0$. As usual, we assume that $0<\varepsilon\ll 1$. 
Our discussion of this example does not provide a new method of solving this problem, but instead it improves the interpretation of the quality of solutions obtained by various methods.

\subsubsection{Regular expansion}

The classical perturbation analysis supposes that the solution to this equation can be written as the power series
\begin{align}
 y(t) = y_0(t) + y_1(t)\e + y_2(t)\e^2+y_3(t)\e^3+\cdots\>.
 \end{align}
Substituting this series in equation \eqref{Duffing} and solving the equations obtained by equating to zero the coefficients of powers of $\e$ in the residual, we find $y_0(t)$ and $y_1(t)$ and we thus have the solution
\begin{align}
z_1(t)= \cos( t) +\e \left( \frac{1}{32}\cos(3t) -\frac{1}{32}\cos( t) -\frac{3}{8}t\sin( t) \right)\>. \label{classical1st}
\end{align}
The difficulty with this solution 
is typically characterized in one of two ways. Physically, the secular term $t\sin t$ shows that our simple perturbative method has failed since the energy conservation prohibits unbounded solutions. Mathematically, the secular term $t\sin t$ shows that our method has failed since the periodicity of the solution contradicts the existence of secular terms.

\begin{figure}
\centering
\includegraphics[width=.55\textwidth]{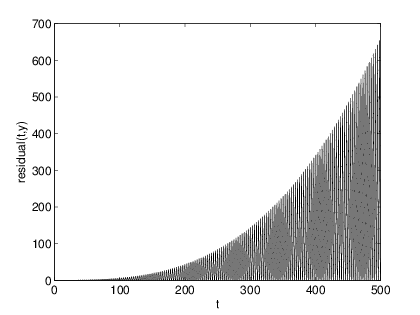}
\caption{Absolute Residual for the first-order classical perturbative solution of the unforced weakly damped Duffing equation with $\e=0.1$.}
\label{ClassicalDuffingRes}
\end{figure}

Both these characterizations are correct, but require foreknowledge of what is physically meaningful or of whether the solutions are bounded. In contrast, interpreting \eqref{classical1st} from the backward error viewpoint is much simpler. To compute the residual, we simply substitute $z_2$ in equation \eqref{Duffing}, that is, the residual is defined by
\begin{align}
\Delta_1(t) = z_1'' + z_1 + \e z_1^3\>.
\end{align}
For the first-order solution of equation \eqref{classical1st}, the residual is
\begin{multline}
\Delta_1(t) =  \Big( -\tfrac {3}{64}\cos( t) +\tfrac{3}{128}\cos( 5t) +\tfrac{3}{128}\cos( 3t) -
\tfrac{9}{32}t\sin(t)\\ -\tfrac{9}{32}t\sin( 3t)\Big) \e^2+O( \e^3) \>. \label{ClassicalRes1st}
\end{multline}
$\Delta_1(t)$ is exactly computable. We don't print it all here because it's too ugly, but in figure \ref{ClassicalDuffingRes}, we see that the complete residual grows rapidly.
 This is due to the secular term $-\tfrac{9}{32}t(\sin(t)-\sin(3t))$ of equation \eqref{ClassicalRes1st}. Thus we come to the conclusion that the secular term contained in the first-order solution obtained in equation \eqref{classical1st} invalidate it, but this time we do not need to know in advance what to physically expect or to prove that the solution is bounded. This is a slight but sometimes useful gain in simplicity.\footnote{In addition, this method makes it easy to find mistakes of various kinds. For instance, a typo in the 1978 edition of \cite{Bender(1978)} was uncovered by computing the residual. That typo does not seem to be in the later editions, so it's likely that the authors found and fixed it themselves.}

A simple Maple code makes it possible to easily obtain higher-order solutions:
\begin{lstlisting}
#Regular Expansion for Duffing's Equation
#We choose initial conditions y(0)=1 and y'(0)=0 so y(t)=cos(t) to O(e).
restart;
macro(e=epsilon);
N := 3;
Order := N+1;
z := add(y[k](t)*epsilon^k, k = 0 .. N);
DE := y -> diff(y, t, t)+y+e*y^3;
des := series( DE(z), e);
dos := dsolve({coeff(des, e, 0), y[0](0) = 1, (D(y[0]))(0) = 0}, y[0](t));
assign(dos);
for k to N do
	tmp:=dsolve({coeff(des,e,k),y[k](0)=0,(D(y[k]))(0)=0},y[k](t));
	assign(tmp);
end do:
Delta := DE(z):
ResidualSeries := map(combine, series(Delta, e, Order+3), trig);
\end{lstlisting}
Experiments with this code suggests the conjecture that $\Delta_n=O(t^n\e^{n+1})$. For this to be small, we must have $\e t=o(1)$ or $t<O(\sfrac{1}{\e})$.

\subsubsection{Lindstedt's method}

\begin{figure}
\centering
\subfigure[First-Order\label{FstLindstedt}]{\includegraphics[width=.48\textwidth]{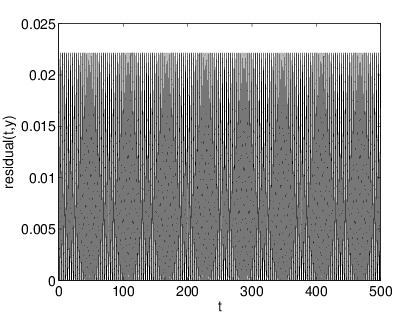}}
\subfigure[Second-Order\label{SndLindstedt}]{\includegraphics[width=.48\textwidth]{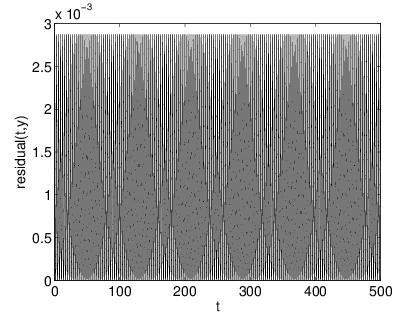}}
\caption{Absolute Residual for the Lindstedt solutions of the unforced weakly damped Duffing equation with $\e=0.1$.}
\end{figure}

The failure to obtain an accurate solution on unbounded time intervals by means of the classical perturbation method suggests that another method that eliminates the secular terms will be preferable. A natural choice is Lindstedt's method, which rescales the time variable $t$ in order to cancel the secular terms. 
The idea is that if we use a rescaling $\tau=\omega t$ of the time variable and chose $\omega$ wisely the secular terms from the classical perturbation method will cancel each other out.\footnote{Interpret this as: we choose $\omega$ to keep the residual small over as long a time-interval as possible.} Applying this transformation, equation \eqref{Duffing} becomes
\begin{align}
\omega^2 y''(\tau)+y(\tau)+\e y^3(\tau) \qquad y(0)=1,\enskip y'(0)=0\>.\label{DuffingTau}
\end{align}
  In addition to writing the solution as a truncated series
\begin{align}
 z_1(\tau) = y_0(\tau)+y_1(\tau)\e \label{ytau}
 \end{align}
we expand the scaling factor as a truncated power series in \e:
\begin{align}
 \omega=1+\omega_1\e\>. \label{omeg}
\end{align}
Substituting \eqref{ytau} and \eqref{omeg} back in equation \eqref{DuffingTau} to obtain the residual and setting the terms of the residual to zero in sequence, we find the equations
\begin{align}
y_0'' + y_0 =0\>,
\end{align}
so that $y_0=\cos(\tau)$, and
\begin{align}
y_1'' + y_1 = -y_0^3 - 2\omega_1 y_0''
\end{align}
subject to the same initial conditions, $y_0(0)=1, y'_0(0)=0, y_1(0)=0$, and $y_1'(0)=0$. By solving this last equation, we find
\begin{align}
y_1(\tau) =\frac{31}{32}\cos(\tau) +\frac{1}{32}\cos(3\tau) -\frac{3}{8}\tau \sin(\tau)+\omega_1\tau\sin(\tau)\>.
\end{align}
So, we only need to choose $\omega_1=\sfrac{3}{8}$ to cancel out the secular terms containing $\tau\sin(\tau)$. Finally, we simply write the solution $y(t)$ by taking the first two terms of $y(\tau)$ and plug in $\tau=(1+\sfrac{3\e}{8})t$:
\begin{align}
z_1(t) = \cos \tau +\e \left( \frac{31}{32}\cos \tau +\frac{1}{32}\cos \tau  \right) 
\end{align}
This truncated power series can be substituted back in the left-hand side of equation \eqref{Duffing} to obtain an expression for the residual:
\begin{align}
\Delta_1(t) =  \left( \frac{171}{128}\cos \left( t \right) +\frac {3}{128}\cos \left( 5t \right) +\frac {9}{16}\cos \left( 3t \right)  \right) \e^2+O \left( \e^3\right) 
\end{align}
See figure \ref{FstLindstedt}.
We then do the same with the second term $\omega_2$. The following Maple code has been tested up to order $12$:
\begin{lstlisting}
# Elimination of Secular Terms in the Solution of the Duffing Equation
# with the Poincare-Lindstedt method.
restart;
macro(e=epsilon);
N := 12;
Order := N+1;
z := add(y[k](tau)*e^k, k = 0..N);
omega := 1+add(a[k]*e^k, k = 1..N);
DE := y -> omega^2*(diff(y, tau, tau))+y+e*y^3;
des := series(DE(z), e);
dos := dsolve({coeff(des,e,0),y[0](0)=1,(D(y[0]))(0)=0},y[0](tau));
assign(dos);
for k to N do
	tmp := convert(combine(coeff(des, e, k), trig), exp);
	UZ := eval(tmp, [exp(I*tau) = Z, exp(-I*tau) = 1/Z]);
	ah := coeff(UZ, Z, 1);
	antisecular := solve(ah = 0, a[k]);
	if {antisecular} <> {} then
		a[k] := antisecular;
	end if;
	tmp := dsolve({evalc(tmp),y[k](0)=0,(D(y[k]))(0)=0},y[k](tau));
	assign(tmp);
end do;
Delta := DE(z);
Sdelta := map(simplify, series(Delta, epsilon, Order+4));
map(combine, Sdelta, trig);
\end{lstlisting}
The significance of this is as follows: The normal presentation of the method first requires a proof (an independent proof) that the reference solution is bounded and therefore the secular term $\e t \sin t$ in the classical solution is spurious. \textsl{But} the residual analysis needs no such proof. It says directly that the classical solution solves not
\begin{align}
 f(t,y,y',y'')=0
 \end{align}
nor $f+\Delta f=0$ for uniformly small $\Delta$ but rather that the residual \textsl{departs} from 0 and is \textsl{not} uniformly small whereas the residual for the Lindstedt solution \textsl{is} uniformly small.

\subsection{Morrison's counterexample}
In \cite[pp.~192-193]{Omalley(2014)}, we find a discussion of the equation
\begin{align}
y''+y+\e(y')^3+3\e^2(y')=0\>.
\end{align}
O'Malley attributed the equation to \cite{Morrison(1966)}. The equation is one that is supposed to illustrate a difficulty with the (very popular and effective) method of multiple scales. We give a relatively full treatment here because a residual-based approach shows that the method of multiple scales, applied somewhat artfully, can be quite successful and moreover we can demonstrate \textsl{a posteriori} that the method was successful. The solution sketched in \cite{Omalley(2014)} uses the complex exponential format, which one of us used to good effect in his PhD, but in this case the real trigonometric form leads to slightly simpler formul\ae. We are very much indebted to our colleague, Professor Pei Yu at Western, for his careful solution, which we follow and analyze here.\footnote{We had asked him to solve this problem using one of his many computer algebra programs; instead, he presented us with an elegant handwritten solution.} 

The first thing to note is that we will use three time scales, $T_0=t$, $T_1=\e t$, and $T_2=\e^2 t$ because the DE contains an $\e^2$ term, which will prove to be important. Then the multiple scales formalism gives
\begin{align}
\frac{d}{dt} = \frac{\partial}{\partial T_0} + \e \frac{\partial}{\partial T_1} + \e^2 \frac{\partial}{\partial T_2} \label{msformalism}
\end{align}
This formalism gives most students some pause, at first: replace an ordinary derivative by a sum of partial derivatives using the chain rule? What could this mean? But soon the student, emboldened by success on simple problems, gets used to the idea and eventually the conceptual headaches are forgotten.\footnote{This can be made to make sense, after the fact. We imagine $F(T_1,T_2,T_3)$ describing the problem, and $\sfrac{d}{dt}=\sfrac{\partial F}{\partial T_1}\sfrac{\partial T_1}{\partial t} + \sfrac{\partial F}{\partial T_2}\sfrac{\partial T_2}{\partial t} + \sfrac{\partial F}{\partial T_3}\sfrac{\partial T_3}{\partial t}$ which gives $\sfrac{d}{dt}=\sfrac{\partial F}{\partial T_1}+\e \sfrac{\partial F}{\partial T_2} + \e^2 \sfrac{\partial F}{\partial T_3}$ if $T_1=t, T_2=\e t$ and $T_3=\e^2t$.} But sometimes they return, as with this example.

To proceed, we take
\begin{align}
y=y_0+\e y_1+\e^2 y_2+O(\e^3)
\end{align}
and equate to zero like powers of $\e$ in the residual. The expansion of $\sfrac{d^2 y}{dt^2}$ is straightforward:
\begin{multline}
\left( \frac{\partial}{\partial T_0} + \e\frac{\partial}{\partial T_1}+\e^2\frac{\partial}{\partial T_2}\right)^2(y_0+\e y_1+\e^2 y_2) =\\
 \frac{\partial^2 y_0}{\partial T_0^2} + \e\left(\frac{\partial^2 y_1}{\partial T_0^2}+2\frac{\partial^2 y_0}{\partial T_0\partial T_1}\right)
+\e^2\left(\frac{\partial^2 y_2}{\partial T_0^2}+2\frac{\partial^2 y_1}{\partial T_0\partial T_1}+\frac{\partial^2 y_0}{\partial T_1^2}+2\frac{\partial^2 y_0}{\partial T_0\partial T_1}\right)
\end{multline}
For completeness we include the other necessary terms, even though this construction may be familiar to the reader. We have
\begin{multline}
\e\left(\frac{dy}{dt}\right)^3 = \e\left( \left(\frac{\partial}{\partial T_0}+\e\frac{\partial}{\partial T_1}\right)(y_0+\e y_1)\right)^3\\
= \e\left(\frac{\partial y_0}{\partial T_0}\right)^3 + 3\e^2\left(\frac{\partial y_0}{\partial T_0}\right)^2 \left(\frac{\partial y_0}{\partial T_1}+\frac{\partial y_1}{\partial T_0}\right)+\cdots\>,
\end{multline}
and $y=y_0+\e y_1+\e^2 y_2$ is straightforward, and also
\begin{align}
3\e^2\left( \left(\frac{\partial}{\partial T_0}+\cdots\right)(y_0+\cdots)\right) = 3\e^2\frac{\partial y_0}{\partial T_0}+\cdots
\end{align}
is at this order likewise straightforward. At $O(\e^0)$ the residual is
\begin{align}
\frac{\partial^2 y_0}{\partial T_0^2}+y_0=0
\end{align}
and without loss of generality we take as solution
\begin{align}
y_0 = a(T_1,T_2)\cos(T_0+\varphi(T_1,T_2))
\end{align}
by shifting the origin to a local maximum when $T_0=0$. For notational simplicity put $\theta=T_0+\varphi(T_1,T_2)$. At $O(\e^1)$ the equation is
\begin{align}
\frac{\partial^2 y_1}{\partial T_0^2} + y_1 = -\left(\frac{\partial y_0}{\partial T_0}\right)^3 - 2\frac{\partial^2 y_0}{\partial T_0\partial T_1}
\end{align}
where the first term on the right comes from the $\e\dot{y}^3$ term whilst the second comes from the multiple scales formalism. Using $\sin^3\theta=\sfrac{3}{4}\sin\theta-\sfrac{1}{4}\sin 3\theta$, this gives
\begin{align}
\frac{\partial^2 y_1}{\partial T_0^2}+y_1 = \left(2\frac{\partial a}{\partial T_1}+\frac{3}{4} a^3\right)\sin\theta + 2a\frac{\partial \varphi}{\partial T_1}\cos\theta - \frac{a^3}{4}\sin 3\theta
\end{align}
and to suppress the resonance that would generate secular terms we put
\begin{align}
\frac{\partial a}{\partial T_1} = -\frac{3}{8}a^3 \quad\textrm{and}\qquad  \frac{\partial\varphi}{\partial T_1}=0\>. \label{525}
\end{align}
Then $y_1 = \frac{a^3}{32}\sin 3\theta$ solves this equation and has $y_1(0)=0$, which does not disturb the initial condition $y_0(0)=a_0$, although since $\sfrac{dy_1}{dT_0}=\sfrac{3a^2}{32}\cos3\theta$ the derivative of $y_0+\e y_1$ will differ by $O(\e)$ from zero at $T_0=0$. This does not matter and we may adjust this by choice of initial conditions for $\varphi$, later.

The $O(\e^2)$ term is somewhat finicky, being
\begin{multline}
\frac{\partial^2 y_2}{\partial T_0^2}+y_2 = -2\frac{\partial^2 y_0}{\partial T_0\partial T_2} -2\frac{\partial^2 y_1}{\partial T_0\partial T_1} \\
 -3\left(\frac{\partial y_0}{\partial T_0}\right)^2 \left(\frac{\partial y_0}{\partial T_1}+\frac{\partial y_1}{\partial T_0}\right) - \frac{\partial^2 y_0}{\partial T_1^2}-3\frac{\partial y_0}{\partial T_0}
\end{multline}
where the last term came from $3(\dot{y})\e^2$. Proceeding as before, and using $\partial\varphi/\partial T_1=0$ and $\sfrac{\partial a}{\partial T_1}=-\sfrac{3}{8}\>a^3$ as well as some other trigonometric identities, we find the right-hand side can be written as
\begin{align}
\left(2\frac{\partial a}{\partial T_2}+3a\right)\sin\theta+\left(2a\frac{\partial\varphi}{\partial T_2}-\frac{9}{128}a^5\right)\cos\theta-\frac{27}{1024}a^5\cos3\theta+\frac{9}{128}a^5\cos5\theta\>.
\end{align}
 Again setting the coefficients of $\sin\theta$ and $\cos\theta$ to zero to prevent resonance we have
 \begin{align}
 \frac{\partial a}{\partial T_2}=-\frac{3}{2}a \label{528}
 \end{align}
 and
 \begin{align}
 \frac{\partial \varphi}{\partial T_2} = \frac{9}{256}a^4\qquad (a\neq0).
 \end{align}
 This leaves
 \begin{align}
 y_2= \frac{27}{1024}a^5\cos3\theta - \frac{3 a^5}{1024}\cos5\theta
 \end{align}
 again setting the homogeneous part to zero.
 
 Now comes a bit of multiple scales magic: instead of solving equations \eqref{525} and \eqref{528} in sequence, as would be usual, we write
 \begin{align}
 \frac{da}{dt} &= \frac{\partial a}{\partial T_0} + \e \frac{\partial a}{\partial T_1} + \e^2\frac{\partial a}{\partial T_2}
= 0 + \e\left(-\frac{3}{8}a^3\right) + \e^2\left(-\frac{3}{2}a\right) \nonumber \\
 &= -\frac{3}{8}\e a(a^2+4\e)\>. \label{magic}
 \end{align}
 Using $a=2R$ this is equation (6.50) in \cite{Omalley(2014)}. Similarly
 \begin{align}
 \frac{d\varphi}{dt} &= \e \frac{\partial\varphi}{\partial T_1}+\e^2 \frac{\partial\varphi}{\partial T_2}
= 0+\e^2 \frac{9}{256} a^4 \label{moremagic}
 \end{align}
 and once $a$ has been identified, $\varphi$ can be found by quadrature. Solving \eqref{magic} and \eqref{moremagic} by Maple,
 \begin{align}
 a = \frac{\sqrt{\e} a_0}{\displaystyle \sqrt{\e e^{3\e^2 t}+\frac{a_0^2}{4}(e^{3\e^2 t}-1)}} = 2\frac{\sqrt{\e} a_0}{\sqrt{u}}
 \end{align}
 and
 \begin{align}
 \varphi = -\frac{3}{16}\e^2\ln u + \frac{9}{16}\e^4 t-\frac{3}{16}\frac{\e^2a_0^2}{u}
 \end{align}
 where $u=4\e e^{3\e^2 t}+a_0^2(e^{3\e^2 t}-1)$. The residual is (again by Maple)
\begin{align}
\footnotesize \e^3\left( \frac{9}{16}a_0^3\cos3t+a_0^7\left( -\frac{351}{4096}\sin t - \frac{9}{512} \sin 7t+\frac{333}{4096}\sin 3t + \frac{459}{4096}\sin 5t\right)\right)+O(\e^4)
\end{align}
and there is no secularity visible in this term.

\begin{figure}
\centering
\includegraphics[width=.55\textwidth]{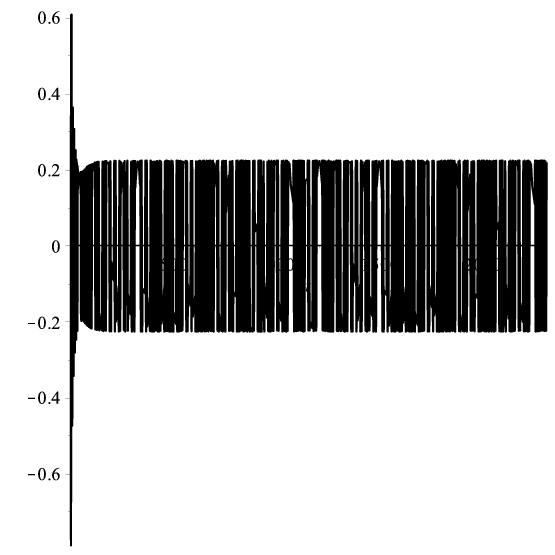}
\caption{The residual $|\Delta_3|$ divided by $\e^3a$, with $\e=0.1$, where $a=O(e^{-\sfrac{3}{2}\>\e^2 t})$, on $0\leq t\leq \sfrac{10\mathrm{ln}(10)}{\e^2}$ (at which point $a=10^{-15}$). We see that $|\sfrac{\Delta_3}{\e^3a}|<1$ on this entire interval.}
\label{YuResidual}
\end{figure}

It is important to note that the construction of the equation \eqref{magic} for $a(t)$ required both $\sfrac{\partial a}{\partial T_1}$ and $\sfrac{\partial a}{\partial T_2}$. Either one alone gives misleading or inconsistent answers. While it may be obvious to an expert that both terms must be used at once, the situation is somewhat unusual and a novice or casual user of perturbation methods may well wish reassurance. (We did!) Computing (and plotting) the residual $\Delta=\ddot{z}+z+\e(\dot{z})^3+3\e^2\dot{z}$ does just that (see figure \ref{YuResidual}).
 It is simple to verify that, say, for $\e=1/100$, $|\Delta|<\e^3a$ on $0<t<10^5\pi$.
Notice that $a\sim O(e^{-\sfrac{3}{2}\>\e^2 t})$ and $e^{-\sfrac{3}{2}\cdot 10^{-4}\cdot 10^5\cdot\pi}=e^{-15\pi} \doteq 10^{-15}$ by the end of this range. The method of multiple scales has thus produced $z$, the exact solution of an equation uniformly and relatively near to the original equation. In trigonometric form,
\begin{multline}
z = a\cos(t+\varphi)+\e\frac{a^3}{32}\cos(3(t+\varphi)) \\
+ \e^2\left(\frac{27}{1024}a^5\cos(3(t+\varphi)) 
 -\frac{3}{1024}a^5\cos^5((5(t+\varphi)) \right) \label{zeqn}
\end{multline} 
and $a$ and $\varphi$ are as in equations \eqref{magic} and \eqref{moremagic}. Note that $\varphi$ asymptotically approaches zero. Note that the trigonometric solution we have demonstrated here to be correct, which was derived for us by our colleague Pei Yu, appears to differ from that given in \cite{Omalley(2014)}, which is
\begin{align}
y= Ae^{it} + \e Be^{3it} + \e^2 Ce^{5it}+\cdots
\end{align}
where (with $\tau=\e t$)
\begin{align}
C\sim \frac{3}{64}A^5+\cdots \qquad \textrm{and}\qquad B\sim -\frac{A^3}{8}(i+\frac{45}{8}\e|A|^2+\cdots)
\end{align}
and, if $A=Re^{i\varphi}$,
\begin{align}
\frac{dR}{d\tau} = -\frac{3}{2}(R^3+\e R+\cdots) 
\qquad \textrm{and}\qquad
\frac{d\varphi}{d\tau} = -\frac{3}{2}R^2 (1+\frac{3\e}{8}R^2+\cdots)
\end{align}
Of course with the trigonometric form $y=a\cos(t+\varphi)$, the equivalent complex form is
\begin{align}
y &= a \left( \frac{e^{it+i\varphi}+ e^{-it-i\varphi}}{2}\right) 
= \frac{a}{2}e^{i\varphi}e^{it}+c.c.
\end{align}
and so $R=\sfrac{a}{2}$. As expected, equation (6.50) in \cite{Omalley(2014)} becomes
\begin{align}
\frac{da}{d\tau}\left(\frac{a}{2}\right) = -\frac{3}{2}\frac{a}{2}\left(\frac{a^2}{4}+\e\right)
\end{align}
or, alternatively,
\begin{align}
\frac{da}{d\tau} = -\frac{3}{8}\e a(a^2+4\e)
\end{align}
which agrees with that computed for us by Pei Yu. However, O'Malley's equation (6.48) gives
\begin{align}
C\cdot e^{i\cdot 5t} &= \frac{3}{64}A^5 e^{i5t} = \frac{3}{64}R^5e^{i5\theta} = \frac{3}{2048}a^5 e^{i5\theta}\>,
\end{align}
so that
\begin{align}
Ce^{i5t}+c.c = \frac{3}{1024}a^5\cos5\theta\>,
\end{align}
whereas Pei Yu has $-\sfrac{3}{1024}$. As demonstrated by the residual in figure \ref{YuResidual}, Pei Yu is correct. Well, sign errors are trivial enough.

More differences occur for $B$, however. The $-\sfrac{A^3}{8} \> ie^{3it}$ term becomes $\sfrac{a^3}{32}\>\cos 3\theta$, as expected, but $-\sfrac{45}{64}A^3\cdot |A|^2e^{3it}+c.c.$ becomes $-\sfrac{45}{32}\sfrac{a^5}{32}\>\cos3\theta = -\sfrac{45}{1024}\>a^5\cos3\theta$, not $\sfrac{27}{1024}\>a^5\cos3\theta$. Thus we believe there has been an arithmetic error in \cite{Omalley(2014)}. This is also present in \cite{Omalley(2010)}.  Similarly, we believe the $\sfrac{d\varphi}{dt}$ equation there is wrong.

Arithmetic errors in perturbation solutions are, obviously, a constant hazard even for experts. We do not point out this error (or the other errors highlighted in this paper) in a spirit of glee---goodness knows we've made our own share. No, the reason we do so is to emphasize the value of a separate, independent check using the residual. Because we have done so here, we are certain that equation \eqref{zeqn} is correct: it produces a residual that is uniformly $O(\e^3)$ for bounded time, and which is $O(\e^{9/2}e^{-\sfrac{3}{2}\>\e^2 t})$ as $t\to \infty$. (We do not know why there is extra accuracy for large times).

Finally, we remark that the difficulty this example presents for the method of multiple scales is that equation \eqref{magic} cannot be solved itself by perturbation methods (or, al least, we couldn't do it). One has to use all three terms at once; the fact that this works is amply demonstrated afterwards.
Indeed the whole multiple scales procedure based on equation \eqref{msformalism} is really very strange when you think about it, but it can be justified afterwards. It really doesn't matter how we find equation \eqref{zeqn}. Once we have done so, verifying that it is the exact solution of a small perturbation of the original equation is quite straightforward. The implementation is described in the following Maple code:
\begin{lstlisting}
restart;
r := epsilon;
de := u -> (diff(u, t, t)+u+r*(diff(u, t))^3+3*r^2*(diff(u, t)));
U := -a0^2+4*exp(3*epsilon^2*t)*epsilon+exp(3*epsilon^2*t)*a0^2;
a := 2*sqrt(r)*a0/sqrt(U);
phi := -(3/16)*epsilon^2*ln(U)+(9/16)*epsilon^4*t-(3/16)*epsilon^2*a0^2/U;
z := a*cos(t+phi)+(1/32)*r*a^3*sin(3*t+3*phi)
+r^2*((27/1024)*a^5*cos(3*(t+phi))-(3/1024)*a^5*cos(5*(t+phi))):
resid := de(z):
zer := MultiSeries[series](resid, r, 4):
map(combine, zer, trig);
eps := 1/10;
plot(eval(resid/(a*r^3), [a0 = 1.0, r = eps]), t = 0 .. 10*ln(10)/eps^2,
 colour = BLACK);
\end{lstlisting}

\subsection{The lengthening pendulum}
As an interesting example with a genuine secular term, \cite{Boas(1966)} discuss the lengthening pendulum. There, Boas solves the linearized equation exactly in terms of Bessel functions. We use the model here as an example of a perturbation solution in a physical context. The original Lagrangian leads to
\begin{align}
\frac{d}{dt} \left(m\ell^2\frac{d\theta}{dt}\right)+mg\ell\sin\theta =0
\end{align}
(having already neglected any system damping). The length of the pendulum at time $t$ is modelled as $\ell =\ell_0+vt$, and implicitly $v$ is small compared to the oscillatory speed $\sfrac{d\theta}{dt}$ (else why would it be a pendulum at all?). The presence of $\sin\theta$ makes this a nonlinear problem; when $v=0$ there is an analytic solution using elliptic functions \cite[chap.~4]{Lawden(2013)}.

We \textsl{could} do a perturbation solution about that analytic solution; indeed there is computer algebra code to do so automatically \cite{Rand(2012)}. For the purpose of this illustration, however, we make the same small-amplitude linerization that Boas did and replace $\sin\theta$ by $\theta$. Dividing the resulting equation by $\ell_0$, putting $\e=\sfrac{v}{\ell_0\omega}$ with $\omega=\sqrt{\sfrac{g}{\ell_0}}$ and rescaling time to $\tau=\omega t$, we get
\begin{align}
(1+\e\tau)\frac{d^2\theta}{d\tau^2}+2\e\frac{d\theta}{d\tau}+\theta=0\>.
\end{align}
This supposes, of course, that the pin holding the base of the pendulum is held perfectly still (and is frictionless besides).

Computing a regular perturbation approximation
\begin{align}
z_{\textrm{reg}} = \sum_{k=0}^N \theta_k(\tau)\e^k
\end{align}
is straightforward, for any reasonable $N$, by using computer algebra. For instance, with $N=1$ we have
\begin{align}
z_{\textrm{reg}} = \cos\tau + \e\left(\frac{3}{4}\sin\tau+\frac{\tau^2}{4}\sin\tau-\frac{3}{4}\tau\cos\tau\right)\>.
\end{align}
This has residual
\begin{align}
\Delta_{\textrm{reg}} &= (1+\e\tau)z''_{\textrm{reg}}+2\e z'_{\textrm{reg}}+z_{\textrm{reg}}\\
&= -\frac{\e^2}{4}\left(\tau^3\sin\tau-9\tau^2\cos\tau-15\tau\sin\tau\right)
\end{align}
also computed straightforwardly with computer algebra. By experiment with various $N$ we find that the residuals are always of $O(\e^{N+1})$ but contain powers of $\tau$, as high as $\tau^{2N+1}$. This naturally raises the question of just when this can be considered ``small.'' We thus have the \textsl{exact} solution of
\begin{align}
(1+\e\tau)\frac{d^2\theta}{d\tau^2}+2\e\frac{d\theta}{d\tau}+\theta = \Delta_{\textrm{reg}}(\tau)=P(\e^{N+1}\tau^{2N+1})
\end{align}
and it seems clear that if $\e^{N+1}\tau^{2N+1}$ is to be considered small it should at least be smaller than $\e\tau$, which appear on the left hand side of the equation. [$\sfrac{d^2}{d\tau^2}$ is $-\cos\tau$ to leading order, so this is periodically $O(1)$.] This means $\e^N\tau^{2N}$ should be smaller than $1$, which forces $\tau\leq T$ where $T=O(\e^{-q})$ with $q<\frac{1}{2}$. That is, this regular perturbation solution is valid only on a limited range of $\tau$, namely, $\tau=O(\e^{-\sfrac{1}{2}})$.

Of course, the original equation contains a term $\e\tau$, and this itself is small only if $\tau\leq T_{\max}$ with $T_{\max}=O(\e^{-1+\delta})$ for $\delta>0$. Notice that we have discovered this limitation of the regular perturbation solution without reference to the `exact' Bessel function solution of this linearized equation. Notice also that $\Delta_{\textrm{reg}}$ can be interpreted as a small forcing term; a vibration of the pin holding the pendulum, say. Knowing that, say, such physical vibrations, perhaps caused by trucks driving past the laboratory holding the pendulum, are bounded in size by a certain amount, can help to decide what $N$ to take, and over which $\tau$-interval the resulting solution is valid.

Of course, one might be interested in the forward error $\theta-z_{\textrm{reg}}$; but then one should be interested in the forward errors caused by neglecting physical vibrations (e.g. of trucks passing by) and the same theory---what a numerical analyst calls a condition number---can be used for both.

But before we pursue that farther, let us first try to improve the perturbation solution. The method of multiple scales, or equivalent but easier in this case the renormalization group method \cite{Kirkinis(2012)} which consists for a linear problem of taking the regular perturbation solution and replacing $\cos\tau$ by $\sfrac{(e^{i\tau}+e^{-i\tau})}{2}$ and $\sin\tau$ by $\sfrac{(e^{i\tau}-e^{-i\tau})}{2i}$, gathering up the result and writing it as $\sfrac{1}{2}\>A(\tau;e)e^{i\tau}+\sfrac{1}{2}\>\bar{A}(\tau;\e)e^{-i\tau}$. One then writes $A(\tau;\e) = e^{L(\tau;\e)}+O(\e^{N+1})$ (that is, taking the logarithm of the $\e$-series for $A(\tau;\e)=A_0(\tau)+\e A_1(\tau)+\cdots+\e^NA_N(\tau)+O(\e^{N+1})$, a straightforward exercise (especially in a computer algebra system) and then (if one likes) rewriting $\sfrac{1}{2}\>e^{L(\tau;\e)+i\tau}+$ c.c. in real trigonometric form again, gives an excellent result. If $N=1$, we get
\begin{align}
\tilde{z}_{\textrm{renorm}}=e^{-\sfrac{3}{4}\>\e\tau}\cos\left(\frac{3}{4}\e+\tau-\e\frac{\tau^2}{4}\right)
\end{align}
which contains an irrelevant phase change $\frac{3}{4}\e$ which we remove here as a distraction to get
\begin{align}
z_{\textrm{renorm}} = e^{-\sfrac{3}{4}\>\e\tau}\cos\left(\tau-\e\frac{\tau^2}{4}\right)\>.
\end{align}
This has residual:
\begin{align}
\Delta_{\textrm{renorm}} &= (1+\e\tau)\frac{d^2z_{\textrm{renorm}}}{d\tau^2} +2\e\frac{dz_{\textrm{renorm}}}{d\tau}+z_{\textrm{renorm}} \nonumber\\
 &= \e^2e^{-\frac{3}{4}\e\tau} \left( (\frac{3}{4}\tau^2-\frac{15}{16})\cos(\tau-\e\frac{\tau^2}{4})-\frac{9}{4}\tau\sin(\tau-\e\frac{t^2}{4})\right)+O(\e^3\tau^3e^{-\frac{3}{4}\e\tau})
\>.
\end{align}
By inspection, we see that this is superior in several ways to the residual from the regular perturbation method. First, it contains the damping term $e^{-\sfrac{3}{4}\>\e\tau}$ just as the computed solution does; this residual will be small compared even to the decaying solution. Second,  at order $N$ it contains only $\tau^{N+1}$ as its highest power of $\e$, not $\tau^{2N+1}$. This will be small compared to $\e\tau$ for times $\tau< T$ with $T=O(\e^{-1+\delta})$ for \emph{any} $\delta>0$; that is, this perturbation solution will provide a good solution so long as its fundamental assumption, that the $\e\tau$ term in the original equation, can be considered `small', is good.

\begin{figure}
\includegraphics[width=.45\textwidth]{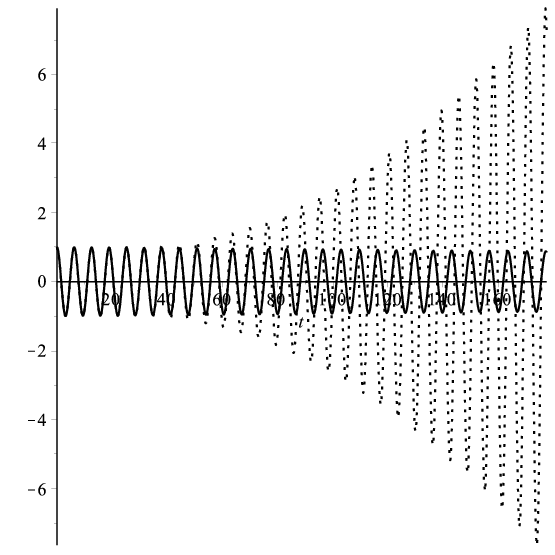}\quad
\includegraphics[width=.45\textwidth]{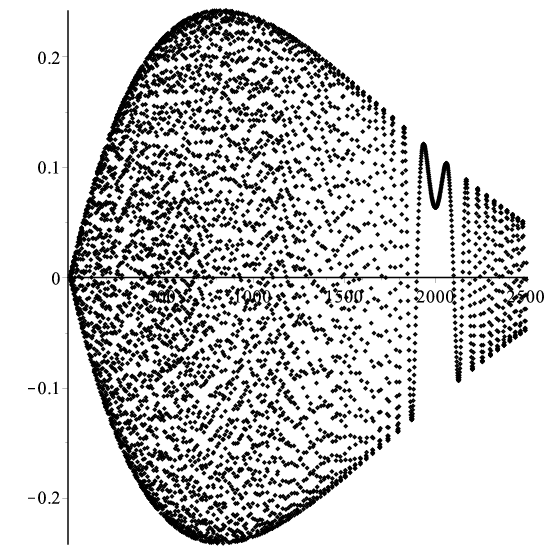}
\caption{On the left, solutions to the lengthening pendulum equation (the renormalized solution is the solid line). On the right, residual of the renormalized solution, which is orders of magnitudes smaller than that of the regular expansion.}
\label{pendulum}
\end{figure}

Note that again the quality of this perturbation solution has been judged without reference to the exact solution, and quite independently of whatever assumptions are usually made to argue for multiple scales solutions (such as boundedness of $\theta$) or the renormalization group method.
Thus, we conclude that the renormalization group method  gives a superior solution in this case, and this judgement was made possible by computing the residual. We have used the following Maple implementation:
\begin{lstlisting}
restart;
macro(e = epsilon);
de := y -> (1+e*t)*(diff(y, t, t))+2*e*(diff(y, t))+y;
z := cos(t);
N := 1;
Order := N+1;
for i to N do
	zt := z+e^i*y[i](t);
	res := series(de(zt), e, i+1);
	eqs := coeff(res, e, i);
	yi := dsolve({eqs, y[i](0) = 0, (D(y[i]))(0) = 0}, y[i](t));
	z := eval(zt, yi);
end do:
res := de(z);
expform := convert(z, exp);
expform := collect(expform, [exp(I*t), exp(-I*t)], factor);
zp := coeff(expform, exp(I*t));
lg := convert(series(ln(series(zp+O(e^Order), e)), e), polynom);
lg := collect(lg, e, factor);
zrg := exp(lg)*exp(I*t);
zrg := zrg+evalc(conjugate(zrg));
zrg := combine(evalc(zrg), trig);
zrg := simplify(zrg);
zrg := exp(-(3/4)*epsilon*t)*cos(t-(1/4)*epsilon*t^2);
resrg := collect(de(zrg), e, t -> combine(simplify(t), trig));
tiny := 1/1000;
plot(eval([z, zrg], e = tiny), t = 0 .. 1/tiny^(3/4), colour = BLACK, linestyle = [2, 1]);
plot(eval(res, e = tiny), t = 1 .. 2500, colour = BLACK, linestyle = 2);
plot(eval(resrg/(tiny*t), e = tiny), t = 1 .. 2500, colour = BLACK,
 style = POINT,numpoints=2016,symbolsize=1);
\end{lstlisting}
%
%
See figure \ref{pendulum}. 

Note that this renormalized residual contains terms of the form $(\e\tau)^k e^{-\sfrac{3}{4}\>\e \tau}$> No matter what order we compute to, these have maxima $O(1)$ when $\tau=O(\sfrac{1}{\e})$, but as noted previously the fundamental assumption of perturbation has been violated by that large a $\tau$.

\paragraph{Optimal backward error again} Now, one further refinement is possible. We may look for an $O(\e^2)$ perturbation of the lengthening of the pendulum, that explains part of this computed residual! That is, we look for $p(t)$, say, so that
\begin{align}
\Delta_2 := (1+\e\tau+\e p(\tau)) z_{\textrm{renorm}}'' + 2(\e+\e^2 p'(\tau))z_{\textrm{renorm}}'+z_{\textrm{renorm}} \label{renormeqs}
\end{align}
has only \textsl{smaller} terms in it  than $\Delta_{\textrm{renorm}}$. Note the correlated changes, $\e^2 p(\tau)$ and $\e^2 p'(\tau)$.

At this point, we don't know if this is possible or useful, but it's a good thing to try. In numerical  analysis terms, we are trying to find a structured backward error for this computed solution.

The procedure for identifying $p(\tau)$ in equation \eqref{renormeqs} is straightforward. We put $p(\tau)=a_0+a_1\tau+a_2\tau^2$ with unknown coefficients, compute $\Delta_2$, and try to choose $a_0$, $a_1$, and $a_2$ in order to make as many coefficients of powers of $\e$ in $\Delta_2$ to be zero as we can. When we do this, we find that
\begin{align}
 p = -\frac{15}{16}+\frac{3}{4}\tau^2
\end{align}
makes
\begin{align}
\Delta_{\textrm{mod}} &= \left(1+\e\tau+\e^2\left(\frac{3}{4}\tau^2-\frac{15}{16}\right)\right)z_{\textrm{renorm}}'' + 2\left(\e+\e^2\left(\frac{3}{2}\tau\right)\right) z_{\textrm{renorm}}' + z_{\textrm{renorm}}\\
&= \e^2e^{-\sfrac{3}{4}\>\e\tau}\left(-\frac{3}{4}\tau\sin\left(\tau-\sfrac{1}{4}\>\e \tau^2\right)\right) + O(\e^3\tau^3 e^{-\sfrac{3\e\tau}{4}})\>.
\end{align}
This is $O(\e^2\tau e^{-\sfrac{3\e\tau}{4}})$ instead of $O(\e^2\tau^2 e^{-\sfrac{3\e\tau}{4}})$, and therefore smaller. This \textsl{interprets} the largest term of the original residual, the $O(\e^2\tau^2)$ term, as a perturbation in the lengthening of the pendulum. The gain is one of interpretation; the solution is the same, but the equation it solves exactly is slightly different. For $O(\e^N\tau^N)$ solutions the modifications will probably be similar. 
Now, if $z\doteq\cos\tau$ then $z'\doteq-\sin\tau$; so if we include a damping term
\begin{align}
\left( +\e^2\cdot\frac{3}{8}\cdot\tau \theta'  \right)
\end{align}
in the model, we have
\begin{align}
\left(1+\e\tau+\e^2\left(\frac{3}{4}\tau^2-\frac{15}{16}\right)\right) z_{\textrm{renorm}}'' + 2\left(\e-\e^2\left(\frac{3}{2}\tau\right)+\e^2\frac{3}{8}\tau\right)z_{\textrm{renorm}}' + z_{\textrm{renorm}} \nonumber\\
= O\left(\e^3\tau^3e^{-\sfrac{3}{4}\>\e\tau}\right)
\end{align}
and \textsl{all} of the leading terms of the residual have been ``explained'' in the physical context. 
If the damping term had been negative, we might have rejected it; having it increase with time also isn't very physical (although one might imagine heating effects or some such).

\subsection{Vanishing lag delay DE}
For another example we consider an expansion that ``everybody knows'' can be problematic. We take the DDE
\begin{align}
\dot{y}(t)+ay(t-\e)+b y(t)=0
\end{align}
from \cite[p.~52]{Bellman(1972)} as a simple instance. Expanding $y(t-\e)=y(t)-\dot{y}(t)\e+O(\e^2)$ we get
\begin{align}
(1-a\e)\dot{y}(t) + (b+a)y(t)=0
\end{align}
by ignoring $O(\e^2)$ terms, with solution
\begin{align}
z(t) = \mathrm{exp}(-\frac{b+a}{1-a\e}t)u_0
\end{align}
if a simple initial condition $u(0)=u_0$ is given. Direct computation of the residual shows
\begin{align}
\Delta &= \dot{z} + az(t-\e)+bz(t)\\
&= O(\e^2)z(t)
\end{align}
uniformly for all $t$; in other words, our computed solution $z(t)$ exactly solves
\begin{align}
\dot{y} + ay(t-\e) + (b+O(\e^2))y(t)=0
\end{align}
which is an equation of the same type as the original, with only $O(\e^2)$ perturbed coefficients. The initial history for the DDE should be prescribed on $-\e\leq t<0$ as well as the initial condition, and that's an issue, but often that history is an issue anyway. So, in this case, contrary to the usual vague folklore that Taylor series expansion in the vanishing lag ``can lead to difficulties'', we have a successful solution and we know that it's successful.

We now need to assess the sensitivity of the problem to small changes in $b$, but we all know that has to be done anyway, even if we often ignore it.

Another example of Bellman's on the same page, $\ddot{y}(t)+ay(t-\e)=0$, can be treated in the same manner. Bellman cautions there that seemingly similar approaches can lead to singular perturbation problems, which can indeed lead to difficulties, but even there a residual/backward error analysis can help to navigate those difficulties.

\subsection{Artificial viscosity in a nonlinear wave equation}

Suppose we are trying to understand a particular numerical solution, by the method of lines, of
\begin{align}
u_t + uu_x = 0 \label{waveeq}
\end{align}
with initial condition $u(0,x)=e^{i\pi x}$ on $-1\leq x\leq 1$ and periodic boundary conditions. Suppose that we use the method of modified equations (see, for example, \cite{Griffiths(1986)}, \cite{Warming(1974)}, or \cite[chap~12]{CorlessFillion(2013)}) to find a perturbed equation that the numerical solution more nearly solves. Suppose also that we analyze the same numerical method applied to the divergence form
\begin{align}
u_t + \frac{1}{2}(u^2)_x=0\>. \label{waveeq2}
\end{align}
Finally, suppose that the method in question uses backward differences $f'(x) = \sfrac{(f(x)-f(x-2\e))}{2\e}$ (the factor 2 is for convenience) on an equally-spaced $x$-grid, so $\Delta x=-2\e$.
The method of modified equations gives
\begin{align}
u_t + uu_x -\e(uu_{xx})+O(\e^2)=0
\end{align}
for equation \eqref{waveeq} and 
\begin{align}
u_t+uu_x -\e (u_x^2 + uu_{xx})+O(\e^2) = 0
\end{align}
for equation \eqref{waveeq2}.

The outer solution to each of these equations is just the reference solution to both equations \eqref{waveeq} and \eqref{waveeq2}, namely,
\begin{align}
u = \frac{1}{i\pi t} W(i\pi t e^{i\pi x})
\end{align}
where $W(z)$ is the principal branch of the Lambert $W$ function, which satisfies $W(z) e^{W(z)}=z$. See \cite{Corless(1996)} for more on the Lambert $W$ function. That $u$ is the solution for this initial condition was first noticed by \cite{weideman(2003)}. 
The residuals of these outer solutions are just $-\e uu_{xx}$ and $-\e(u_x^2+uu_{xx})$ respectively. Simplifying, and again suppressing the argument of $W$ for tidiness, we find that
\begin{align}
-\e uu_{xx} = -\frac{\e W^2}{t^2(1+W^3)}
\end{align}
and
\begin{align}
-\e(u_x^2 + uu_{xx}) = -\frac{\e W^2(2+W)}{t^2(1+W^3)}
\end{align}
where $W$ is short for $W(i\pi t e^{i\pi x})$. We see that if $x=\sfrac{1}{2}$ and $t=\sfrac{1}{(\pi e)}$, both of these are singular:
\begin{align}
-\e uu_{xx} \sim -\e \left( \frac{i\pi^2 e^2\sqrt{2}}{4(et\pi-1)^{\sfrac{3}{2}}}+O\left(\frac{1}{et\pi-1}\right)\right)
\end{align}
and
\begin{align}
-\e (u^2_x +uu_{xx}) \sim -\e \left( \frac{i\pi^2e^2\sqrt{2}}{4(et\pi-1)^{\sfrac{3}{2}}} + O\left(\frac{1}{\sqrt{et\pi-1}}\right)\right)\>.
\end{align}
We see that the outer solution makes the residual very large near $x=\sfrac{1}{2}$ as $t\to \sfrac{1}{(\pi e)}^-$ suggesting that the solution of the modified equation---and thus the numerical solution---will depart from the outer solution. Both the original form and the divergence form are predicted to have similar behaviour, and this is confirmed by numerical experiments.

We remark that using forward differences instead just changes the sign of $\e$, and given the similarity of $euu_{xx}$ to $\e u_{xx}$, we intuit that this will blow up rather quickly, like the backward heat equation, because the exact solution to Burger's equation $u_t+uu_x=\e u_{xx}$ involves a change in variable to the heat equation \cite[pp.~352-353]{Kevorkian(2013)}. We also remark also that this use of residual is a bit perverse: we here substitute the reference solution into an approximate (reverse-engineered) equation. Some authors do use `residual' or even `defect' in this sense., e.g., \cite{Chiba(2009)}. It only fits our usage because the reference solution to the original equation is just the outer solution of the perturbation problem of interest here.

Finally, we can interpolate the numerical solution using a trigonometric interpolant in $x$ tensor producted with the interpolant in $t$ provided by the numerical solver (e.g., \texttt{ode15s} in Matlab). We can then compute the residual $\Delta(t,x) = z_t+zz_x$ in the original equation and we find that, away from the singularity, it is $O(\e)$. If we compute the residual in the modified equation
\begin{align}
\Delta_1(t,x)=z_t+zz_x-\e zz_{xx}
\end{align}
we find that, away from the singularity, it is $O(\e^2)$. This is a more traditional use of residual in a numerical computation, and is done without knowledge of any reference solution. The analogous use we are making for perturbation methods can be understood from this numerical perspective.

\section{Concluding Remarks}
Decades ago, van Dyke had already made the point that, in perturbation theory, ``[t]he possibilities are too diverse to be subject to rules'' \cite[p.~31]{vanDyke(1964)}. Van Dyke was talking about the useful freedom to choose expansion variables artfully, but the same  might be said for perturbation methods generally. This paper has attempted (in the face of that observation) to lift a known technique, namely the residual as a backward error, out of numerical analysis and apply it to perturbation theory. The approach is surprisingly useful and clarifies several issues, namely 
\begin{itemize}
\item BEA allows one to directly use approximations taken from divergent series in an optimal fashion without appealing to ``rules of thumb'' such as stopping before including the smallest term.
\item BEA allows the justification of removing spurious secular terms, even when true secular terms are present.
\item Not least, residual computation and \emph{a posteriori} BEA makes detection of slips, blunders, and bugs all but certain, as illustrated in our examples.
\item Finally BEA interprets the computed solution solution $z$ as the exact solution to just as good a model.
\end{itemize}
%
%
%
In this paper we have used BEA to demonstrate the validity of solutions obtained by the iterative method, by Lindstedt's method, by the method of multiple scales, by the renormalization group method, and by matched asymptotic expansions.
We have also successfully used the residual and BEA in many problems not shown here: eigenvalue problems from \cite{Nayfeh(2011)}; an example from \cite{vanDyke(1964)} using the method of strained coordinates; and many more.

The examples here have largely been for algebraic equations and for ODEs, but the method was used to good effect in \cite{Corless(2014)} for a PDE system describing heat transfer  between concentric cylinders, with a high-order perturbation series in Rayleigh number. Aside from the amount of computational work required, there is no theoretical obstacle to using the technique for other PDE; indeed the residual of a computed solution $z$ (perturbation solution, in this paper) to an operator equation $\varphi(y;x)=0$ is usually computable: $\Delta = \varphi(z;x)$ and its size (in our case, leading term in the expansion in the gauge functions) easily assessed.

It's remarkable to us that the notion, while present here and there in the literature, is not used more to justify the validity of the perturbation series.

We end with a caution. Of course, BEA is not a panacea. There are problems for which it is not possible. For instance, there may be hidden constraints, something like solvability conditions, that play a crucial role and where the residual tells you nothing. A residual can even be zero and if there are multiple solutions, one needs a way to get the right one. 
%
%
%
%
%
%
There are things that can go wrong with this backward error approach. First, the final residual computation might not be independent enough from the computation of $z$, and repeat the same error. An example is if one correctly solves
\begin{align}
\ddot{y}+y+\e \dot{y}^3+3\e^2\dot{y}=0
\end{align}
and verifies that the residual is small, while \textsl{intending} to solve
\begin{align}
\ddot{y}+y+\e \dot{y}^3-3\e^2\dot{y}=0\>,
\end{align}
i.e., getting the wrong sign on the $\dot{y}$ term, both times. Another thing that can go wrong is to have an error in your independent check but not your solution. This happened to us with 183 instead of 138 in subsection \ref{systems}; the discrepancy alerted us that there \textsl{was} a problem, so this at least was noticeable. A third thing that can go wrong is that you verify the residual is small but forget to check the boundary counditions. A fourth thing that can go wrong is that the residual may be small in an absolute sense but still larger than important terms in the equation---the residual may need to be smaller than you expect, in order to get good qualitative results. A fifth thing is that the residual may be small but of the `wrong character', i.e., be unphysical. Perhaps the method has introduced the equivalent of negative damping, for instance. This point can be very subtle.

A final point is that a good solution needs not just a small backward error, but also information about the sensitivity (or robustness) of the model to physical perturbations. We have not discussed computation of sensitivity, but we emphasize that even if $\Delta\equiv 0$, you still have to do it, because real situations have real perturbations. Nonetheless, we hope that we have convinced you that BEA can be helpful.

\bibliographystyle{plain}
\bibliography{perturbationbeaarxiv}

\begin{thebibliography}{10}

\bibitem{Avrachenkov(2013)}
Konstantin~E. Avrachenkov, Jerzy~A. Filar, and Phil~G. Howlett.
\newblock {\em Analytic perturbation theory and its applications}.
\newblock SIAM, 2013.

\bibitem{Bellman(1972)}
Richard~E. Bellman.
\newblock {\em Perturbation techniques in mathematics, physics, and
  engineering}.
\newblock Dover Publications, 1972.

\bibitem{Bender(1978)}
C.M. Bender and S.A. Orszag.
\newblock {\em Advanced mathematical methods for scientists and engineers:
  Asymptotic methods and perturbation theory}, volume~1.
\newblock Springer Verlag, 1978.

\bibitem{Boas(1966)}
Mary~L. Boas.
\newblock {\em Mathematical Methods in the Physical Sciences}.
\newblock John Wiley, New York, 1966.

\bibitem{Boyd(2014)}
John~P. Boyd.
\newblock {\em Solving Transcendental Equations}.
\newblock SIAM, 2014.

\bibitem{Chiba(2009)}
Hayato Chiba.
\newblock Extension and unification of singular perturbation methods for odes
  based on the renormalization group method.
\newblock {\em SIAM Journal on Applied Dynamical Systems}, 8(3):1066--1115,
  2009.

\bibitem{Corless(1993)b}
Robert~M. Corless.
\newblock What is a solution of an {ODE}?
\newblock {\em ACM SIGSAM Bulletin}, 27(4):15--19, 1993.

\bibitem{Corless(1992)}
Robert~M. Corless and George~F. Corliss.
\newblock Rationale for guaranteed {ODE} defect control.
\newblock In L.~Atanassova and J.~Herzberger, editors, {\em Computer Arithmetic
  and Enclosure Methods}, pages 3--12. North-Holland, 1992.

\bibitem{CorlessFillion(2013)}
Robert~M. Corless and N.~Fillion.
\newblock {\em A Graduate Introduction to Numerical Methods, From the Viewpoint
  of Backward Error Analysis}.
\newblock Springer, New York, 2013.
\newblock 868pp.

\bibitem{Corless(1996)}
Robert~M. Corless, G.H. Gonnet, D.E.G. Hare, D.J. Jeffrey, and Donald~E. Knuth.
\newblock On the {Lambert $W$} function.
\newblock {\em Advances in Computational Mathematics}, 5(1):329--359, 1996.

\bibitem{Bruijn(1981)}
Nicolaas~Govert De~Bruijn.
\newblock {\em Asymptotic methods in analysis}, volume~4.
\newblock Dover, 1970.

\bibitem{Deuflhard(2003)}
P.~Deuflhard and A.~Hohmann.
\newblock {\em Numerical analysis in modern scientific computing: an
  introduction}, volume~43.
\newblock Springer Verlag, 2003.

\bibitem{NIST:DLMF}
{NIST Digital Library of Mathematical Functions}.
\newblock http://dlmf.nist.gov/, Release 1.0.10 of 2015-08-07.

\bibitem{Enright(1989)b}
Wayne~H. Enright.
\newblock Analysis of error control strategies for continuous runge-kutta
  methods.
\newblock {\em SIAM Journal on Numerical Analysis}, 26(3):588--599, 1989.

\bibitem{Enright(1989)a}
Wayne~H. Enright.
\newblock A new error-control for initial value solvers.
\newblock {\em Applied Mathematics and Computation}, 31:288--301, 1989.

\bibitem{Geddes(1992)b}
K.~O. Geddes, S.~R. Czapor, and G.~Labahn.
\newblock {\em Algorithms for computer algebra}.
\newblock Kluwer Academic, Boston, 1992.

\bibitem{Grcar(2011)}
J.F. Grcar.
\newblock John von {N}eumann's analysis of {G}aussian elimination and the
  origins of modern numerical analysis.
\newblock {\em SIAM review}, 53(4):607--682, 2011.

\bibitem{Griffiths(1986)}
D.F. Griffiths and J.-M. Sanz-Serna.
\newblock On the scope of the method of modified equations.
\newblock {\em SIAM Journal on Scientific and Statistical Computing},
  7(3):994--1008, 1986.

\bibitem{Higham(1996)}
Nicholas~J. Higham.
\newblock {\em Accuracy and Stability of Numerical Algorithms}.
\newblock SIAM, Philadelphia, 2nd edition, 2002.

\bibitem{Holmes(1995)}
M.H. Holmes.
\newblock {\em Introduction to perturbation methods}.
\newblock Springer, 1995.

\bibitem{Kevorkian(2013)}
Jirair Kevorkian and Julian~D Cole.
\newblock {\em Perturbation methods in applied mathematics}.
\newblock Springer, 2013.

\bibitem{Kirkinis(2012)}
Eleftherios Kirkinis.
\newblock The renormalization group: A perturbation method for the graduate
  curriculum.
\newblock {\em SIAM Review}, 54(2):374--388, 2012.

\bibitem{Lanczos(1988)}
Cornelius Lanczos.
\newblock {\em Applied analysis}.
\newblock Dover Pubns, 1988.

\bibitem{Lawden(2013)}
Derek~F. Lawden.
\newblock {\em Elliptic functions and applications}, volume~80.
\newblock Springer Science \& Business Media, 2013.

\bibitem{Morrison(1966)}
J.A. Morrison.
\newblock Comparison of the modified method of averaging and the two variable
  expansion procedure.
\newblock {\em SIAM Review}, 8(1):66--85, 1966.

\bibitem{Nayfeh(2011)}
Ali~H Nayfeh.
\newblock {\em Introduction to perturbation techniques}.
\newblock John Wiley \& Sons, 2011.

\bibitem{OEIS}
{The On-Line Encyclopedia of Integer Sequences}.
\newblock https://oeis.org/.

\bibitem{Omalley(2014)}
Robert~E. O'Malley.
\newblock {\em Historical Developments in Singular Pertubations}.
\newblock Springer, 2014.

\bibitem{Omalley(2010)}
Robert~E. O'Malley and Eleftherios Kirkinis.
\newblock A combined renormalization group-multiple scale method for singularly
  perturbed problems.
\newblock {\em Studies in Applied Mathematics}, 124(4):383--410, 2010.

\bibitem{Rand(2012)}
Richard Rand and Dieter Armbruster.
\newblock {\em Perturbation methods, bifurcation theory and computer algebra},
  volume~65.
\newblock Springer Science \& Business Media, 2012.

\bibitem{Salvy(2010)}
Bruno Salvy and John Shackell.
\newblock Measured limits and multiseries.
\newblock {\em Journal of the London Mathematical Society}, 82(3):747--762,
  2010.

\bibitem{vanDyke(1964)}
Milton Van~Dyke.
\newblock {\em Perturbation methods in fluid mechanics}.
\newblock Academic Press, 1964.

\bibitem{Warming(1974)}
R.F. Warming and B.J. Hyett.
\newblock The modified equation approach to the stability and accuracy analysis
  of finite-difference methods.
\newblock {\em Journal of computational physics}, 14(2):159--179, 1974.

\bibitem{weideman(2003)}
J.A.C. Weideman.
\newblock Computing the dynamics of complex singularities of nonlinear {PDE}s.
\newblock {\em SIAM J. Appl. Dyn. Syst}, 2(2):171--186, 2003.

\bibitem{Wilkinson(1963)}
James~H. Wilkinson.
\newblock {\em Rounding Errors in Algebraic Processes}.
\newblock Prentice-Hall Series in Automatic Computation. Prentice-Hall,
  Englewood Cliffs, 1963.

\bibitem{Wilkinson(1965)}
James~H. Wilkinson.
\newblock {\em The Algebraic Eigenvalue Problem}.
\newblock Oxford University Press, New York, 1965.

\bibitem{Wilkinson(1971)}
James~H. Wilkinson.
\newblock Modern error analysis.
\newblock {\em SIAM Review}, 13(4):548--568, 1971.

\bibitem{Wilkinson(1984)}
James~H. Wilkinson.
\newblock The perfidious polynomial.
\newblock In Gene~H. Golub, editor, {\em Studies in Numerical Analysis},
  volume~24, pages 1--28. Mathematical Assosication of America, 1984.

\bibitem{Corless(2014)}
Yiming Zhang and Robert~M. Corless.
\newblock High-accuracy series solution for two-dimensional convection in a
  horizontal concentric cylinder.
\newblock {\em SIAM Journal on Applied Mathematics}, 74(3):599--619, 2014.

\end{thebibliography}

\end{document}